\documentclass[12pt]{article}
\usepackage[english,russian]{babel}
\usepackage[cp1251]{inputenc}
\usepackage{amssymb,latexsym,amsmath}
\usepackage{amsfonts}
\usepackage{indentfirst}
\usepackage{graphicx,graphics,verbatim}
\usepackage{cite}
\usepackage{cmap}

\newtheorem{theorem}{Теорема}
\newtheorem{definition}{Определение}
\newtheorem{corollary}{Следствие}

\begin{document}

\begin{center}		
	{\textbf{ On generalizations of discrete and integral Cauchy--Bunyakovskii inequalities
by the method of mean values. Some applications.}}
\end{center}	

\begin{center}
	\textbf{Sergei M. Sitnik}\\
\null
	 Belgorod State National Research     University (BelGU),\\
 Pobedy Street, 85, 308015 Belgorod, Russia\\
	sitnik@bsu.edu.ru\\
\end{center}

\null

{MSC 26D15}

\null
	

\null

{\bf Abstract:\\}

In this preprint we consider generalizations of discrete and integral Cauchy--Bunyakovskii inequalities by the method of mean values with some applications. Mostly the material is compiled as a short survey but some results are proved. Main results are listed, including an interesting inequality with maximum and minimum values. Some applications are considered from different fields of mathematics. Among them are estimates for some special functions, including Euler gamma and incomplete gamma function, the Legendre complete elliptic integrals of the first kind. Also some further possible generalizations are considered and outlined, including generalizations of the Acz\'el and Minkovskii inequalities, a case of spaces with sign--indefinite form, the Jackson's $q$-- integrals.

\newpage

\begin{center}		
\textbf{ Об обобщениях дискретного и интегрального неравенств Коши--Буняковского\\
 методом средних значений.\\ Некоторые приложения.}
\end{center}

\begin{center}
	\textbf{Ситник С.М.}\\
\null
	Белгородский государственный национальный исследовательский университет\\
 (НИУ "БелГУ")\\
 ул. Победы, 85, 308015, Белгород, Россия.\\
	sitnik@bsu.edu.ru\\
\end{center}

\null

{ УДК 517.16,  517.165}

\null

{\bf Ключевые слова: } неравенства Коши--Буняковского, средние значения, неравенство Минковского,  неравенство Ацеля, специальные функции, эллиптический интеграл Лежандра, гамма-- функция, $q$-- интеграл Джексона.

\null

{\bf Аннотация: \\}
В статье в форме небольшого обзора рассматриваются уточнения в терминах средних значений интегрального и дискретного неравенств Коши--Буняковского. Сформулированы основные результаты, в том числе приведено  одно из таких уточнений, содержащее в качестве выбранных средних максимум и минимум функций. Рассматриваются некоторые приложение к получению оценок для специальных функций, а именно, эллиптических интегралов Лежандра и гамма---функций. Также обсуждается  возможность некоторых дальнейших обобщений, в том числе на случай неравенств  Минковского, пространств со знаконеопределённой формой, $q$-- интегралов Джексона.

\newpage

\today\\

\vspace*{3mm}
{ March 23, 2022}

\null

\tableofcontents

\newpage

\section{Введение}

Неравенства Коши--Буняковского для конечных сумм, рядов и интегралов относятся к числу классических, они используются практически во всех разделах математики, как теоретических, так и прикладных. Из ссылок ограничимся указаниями на хорошо известные общие монографии по неравенствам \cite{BB,HLP,MO2,MPF}, а также на прекрасно написанные специальные монографии Севера Драгомира \cite{Dr1} и Джона Майкла Стила \cite{Ste}. Отметим также работы автора \cite{Sit1,Sit3,Sit4,Sit2}.

Неравенство Коши--Буняковского для сумм было доказано О.~Коши в 1821 г., а неравенство Коши--Буняковского для интегралов было опубликовано В.Я.\,Буняковским в 1859 г., см. \cite{Bun},
 и затем переоткрыто через двадцать пять лет Г.А.\,Шварцем в 1884.

 Известно много обобщений и уточнений неравенства Коши--Буняковского, см., например, \cite{BB,Sit1,Sit3,Sit4,HLP,Dr1,MO2,MPF,Sit2,Ste}. В работах \cite{Sit1,Sit3,Sit4,Sit2} изложен принадлежащий автору метод средних значений для построения некоторых классов таких уточнений. Кратко говоря, суть этого метода заключается в том, что по каждому абстрактному среднему, удовлетворяющему набору естественных аксиом, выписываются в явном виде уточнения дискретного или интегрального неравенств Коши--Буняковского. Указанная методика возникла из анализа одного неравенства Милна \cite{HLP} и теоремы Карлица--Дэйкина--Элиезера (теорема CDE) для дискретного случая \cite{Sit1,Dr1,Sit2,Ste}. В работе кратко изложены  основные результаты этого метода, следуя работам автора \cite{Sit1,Sit3,Sit4,Sit2}.

 Указанные неравенства находят многочисленные приложения при изучении дифференциальных уравнений,  в теории операторов преобразования, в том числе  при оценках ядер операторов преобразования,  соответствующих специальных функций, норм операторов преобразования в функциональных пространствах, см., например,  \cite{T1,T2,T4,T3}.

\section{Средние значения}

Различные вопросы теории средних  достаточно подробно изложены в
литературе \cite{BB,Jini,Kal,Sit1,HLP,Fum,MO2,MPF,MBV,Sit2,Toad}. Особенно отметим  специально
посвящённые этим вопросам монографии \cite{Jini,Fum,MBV,Toad}, а среди них выделим фундаментальный труд Коррадо Джини \cite{Jini}, описание происхождения средних в трудах греческих математиков из теории пропорций \cite{Jini,Toad}, а также интересные
приложения средних в теории операторов \cite{Fum}.

\subsection{Аксиоматическое определение абстрактных средних}

\begin{definition}
Пусть даны два числа $x\ge 0, y\ge 0$. Назов\"{е}м функцию $M(x,y)$ {\it средним} этих чисел, если для не\"{е} выполнены следующие свойства (аксиомы).

\begin{enumerate}
\item{Свойство промежуточности:
$$
min(x,y)\le M(x,y) \le max(x,y).
$$}
\item{Свойство несмещ\"{е}нности:
$$
M(x,x)=x.
$$}
\item{Свойство однородности:
$$
M(ax,ay)=a\cdot M(x,y).
$$}
\item{Свойство монотонности: если $x_1 \le x_2$ и $y_1 \le y_2$, то
$$
M(x_1,y) \le M(x_2,y), \ \ \ M(x,y_1) \le M(x,y_2).
$$}
\end{enumerate}
\end{definition}

Аксиоматическое определение средних используется во многих работах, начиная с Коши. Отметим, что свойства симметричности среднего для наших результатов не требуется, допускаются несимметричные средние.

Среди множества средних наиболее известными являются средние арифметическое, геометрическое, квадратичное, полуквадратичное, гармоническое:
\begin{equation*} A(x,y)=\frac{x+y}{2}, G(x,y)=\sqrt{x,y}, Q(x,y)=\sqrt{\frac{x^2+y^2}{2}},
\end{equation*}
\begin{equation*}
S(x,y)=\left(
\frac{\sqrt{x}+\sqrt{y}}{2}\right)^2, H(x,y)=A(\frac{1}{x},\frac{1}{y})=\frac{2xy}{x+y}, \ x,y>0.
\end{equation*}

Нам также понадобится понятие дополнительного среднего.

\begin{definition}
Для среднего $M(x,y)$ {\it дополнительное} к нему среднее определяется по формуле
\begin{equation}
M^*(x,y) = \frac{xy}{M(x,y)};\  x,y>0. \label{1}
\end{equation}
\end{definition}

\subsection{Степенные средние }

 Первым из наиболее известных классов средних, зависящих от параметра, являются степенные средние
\begin{equation}
M(x,y) = M_{\alpha} (x,y) = \left( \frac{x^{\alpha}+y^{\alpha}}{2}
\right)^{\frac{1}{\alpha}} , -\infty < \alpha < \infty \ ,
\alpha \neq 0\ ;\label{2}
\end{equation}
$$
M_{- \infty} (x,y) = \min (x, y)\ ,\  M_{0} = \sqrt{xy}\ , \
M_{\infty} (x, y) = \max (x, y)\ .
$$
Три исключительных значения получаются из неисключительных предельным переходом. Аксиомы абстрактного среднего  проверяются непосредственно.

Для степенных средних от положительных чисел дополнительное среднее находится по формуле
$$
(M_\alpha)^*=M_{-\alpha}.
$$

Степенные средние широко используются в Анализе, в качестве
общеизвестных примеров упомянем только обобщённое суммирование рядов с
использованием средних арифметических (простейшая форма метода
Чезаро) и теорему Фейера о равномерной сходимости средних арифметических частичных сумм ряда
Фурье для непрерывной функции.

\subsection{Средние Тибора Радо}

 Вторым известным классом параметрических средних являются средние Т. Радо
\begin{equation}
R_{\beta}(x,y) = \left( \frac{x^{\beta + 1} - y^{\beta + 1}}{(\beta
+ 1)(x-y)}\right) ^{\frac{1}{\beta}} , -\infty < \beta <
\infty ,\  \beta \neq 0, -1;\label{3}
\end{equation}
$$ R_{- \infty} (x,y) = \min (x, y),  R_{\infty} (x, y) = \max (x, y).
$$
Исключительные значения приводят к логарифмическому среднему
\begin{equation}
R_{-1}(x,y)=L(x,y)=\frac{y-x}{\ln y-\ln x} \label{4}
\end{equation}
и "многоэтажно--показательному" (identric  в англоязычной литературе)
\begin{equation}
R_{0}(x,y)= \frac{1}{e} \left(
\frac{y^{y}}{x^{x}}\right)^{\frac{1}{y-x}}. \label{5}
\end{equation}
Среднее Радо при $\beta=2$ называется средним Герона.

Понятно, что подобные величины встречались намного раньше, но средние $R_{\beta}(x,y)$ были введены как отдельный объект и подробно изучены Тибором Радо в \cite{Rado}, см.  также   \cite{Sit1,Sit3,Sit4,Har,MPF,Sit2}.

 Средние Радо также образуют шкалу по параметру:
$$
\beta_{1} > \beta_{2} \Rightarrow R_{\beta_{1}} (x,y) \geq
R_{\beta_{2}} (x,y),~  \forall x, y.
$$
Четыре исключительных значения  $\beta = \{-\infty, -1, 0,
+\infty\}$ могут быть получены из неисключительных предельным
переходом. Аксиомы абстрактного среднего  проверяются
непосредственно.

Дополнительное среднее считается по простой формуле
$$
(R_{\beta}(x,y)^*=\frac{1}{R_{\beta}(1/x, 1/y)}.
$$

Тибор Радо провел в \cite{Rado} детальное исследование средних $R_\beta$. В
частности‚ им были доказаны замечательные теоремы о связи двух
основных шкал средних  $M_\alpha$ и $R_\beta$.
Первый результат устанавливает, какие средние входят одновременно в обе шкалы: степенных средних и средних Радо. Оказывается, что это только пять основных классических средних.

\begin{theorem} (Т. Радо, 1935). Только следующие пять средних входят одновременно в обе шкалы средних: степенных и Радо, при этом
$$
M_{-\infty} = R_{-\infty}, M_{0} = R_{-2}, M_{\frac{1}{2}} =
R_{\frac{1}{2}}, M_{1} = R_{1}, M_{\infty} = R_{\infty}.
$$
\end{theorem}
Наиболее часто используемые средние ($\min$‚ $\max$‚ арифметическое‚
геометрическое, полуквадратичное) и только они входят в обе шкалы.

Другим замечательным результатом
Радо является полное описание множеств параметров  $(\alpha,\beta)$‚
при которых выполнены неравенства
\begin{equation} \label{est}
M_{\alpha_1}\leq R_\beta \leq M_{\alpha_2}, R_{\beta_1}\leq M_\alpha
\leq R_{\beta_2},
\end{equation}

К сожалению, условия  теорем Радо, при которых выполнены оценки \eqref{est},
приведены в \cite{Rado} (также как и потом в \cite{MPF}) не в явном виде. Расшифровка этих условий не в
терминах некоторых дополнительных неравенств и специально определяемых множеств, а непосредственно
через параметры средних $\alpha, \beta$, получена в  \cite{Sit1,Sit3,Sit4,Sit2}. Результаты содержит

\begin{theorem}\label{Rado2} (Т. Радо, 1935). Справедливы следующие двусторонние неулучшаемые
оценки средних Радо через степенные средние:
$$
M_{\frac{\alpha + 2}{3}}\leq R_{\alpha}\leq M_{0}  \mbox{\ , при \ }
\alpha \in (-\infty, -2],
$$
$$
M_{0}\leq R_{\alpha}\leq M_{\frac{\alpha +2}{3}}  \mbox{\ , при \ }
\alpha \in [-2, -1],
$$
$$
M_{\frac{\alpha\ln2}{\ln(1+\alpha)}}\leq R_{\alpha} \leq
M_{\frac{\alpha +2}{3}} \mbox{\ , при \ }  \alpha \in (-1, -1/2],
$$
$$
M_{\frac{\alpha+2}{3}}\leq R_{\alpha} \leq M_{\frac{\alpha \ln
2}{\ln (1+\alpha)}} \mbox{\ , при \ } \alpha \in [-1/2, 1),\\
$$
\mbox{(при $\alpha = 0$  последнее неравенство понимается в
предельном смысле}\\ $M_{\frac{2}{3}}\leq R_0 \leq M_{\ln 2}$),
$$M_{\frac{\alpha\ln
2}{\ln (1+\alpha)}}\leq R_\alpha \leq M_{\frac{\alpha +2}{3}}
\mbox{\ , при \ } \alpha \in [1,\infty].\\
$$
\end{theorem}

Разумеется, приведённые выше неравенства можно переписать в другом порядке, чтобы они превратились в точные двусторонние оценки степенных средних через средние Радо.

В частном случае $\alpha=-1$‚ из  результатов Т. Радо, приведённых в теореме \ref{Rado2},  следует‚ что фольклорное неравенство для среднего логарифмического
\begin{equation} \label{log1}
M_{0} \leq L \leq M_{1},
\end{equation}
которое  содержится еще в пионерской работе В.Я. Буняковского \cite{Bun}
в качестве примера на приложения его интегрального неравенства,
может быть усилено до следующего:
\begin{equation}\label{log2}
M_{0}(x,y) = \sqrt{xy}\leq L(x,y) = \frac{x-y}{\ln x - \ln y}\leq
M_{\frac{1}{3}}(x,y) = \left(
\frac{x^{\frac{1}{3}}+y^{\frac{1}{3}}}{2}\right)^{3},
\end{equation}
причем порядки средних $0$ и $1/3$ являются неулучшаемыми. Последнее
неравенство многократно переоткрывалось в литературе, а
также приводилось на различных олимпиадах без указания авторства. Не менее эффектно выглядят и неравенства Радо для многоэтажно--показательного среднего
\begin{eqnarray*}
M_{\frac{2}{3}}(x,y) = \left(
\frac{x^{\frac{2}{3}}+y^{\frac{2}{3}}}{2}\right)^{\frac{2}{3}} \leq
R_{0}(x,y) = \\ = \frac{1}{e}\left(
\frac{y^{y}}{x^{x}}\right)^{\frac{1}{y-x}} \leq M_{\ln2}(x,y) =
\left( \frac{{x^{\ln2}}+y^{\ln2}}{2}\right)^{\frac{1}{\ln2}},
\end{eqnarray*}
порядки средних в которых тоже неулучшаемы (простейшая оценка
$M_0\leq R_0 \leq M_1$ также содержится в указанной работе
Буняковского \cite{Bun}). \ Отметим‚
что дальнейшее развитие результатов Радо для выпуклых и вогнутых
функций получено в \cite{Har}. Отметим также, что из известного  неравенства Эрмита - Адамара для экспоненциальной
функции следует только простейшая оценка \eqref{log1}
 $M_0 \leq L \leq M_1$, что  хуже точных неравенств Радо \eqref{log2}.

Средние Радо имеют простой аналитический смысл --- это промежуточные
значения в теореме Лагранжа о средних для логарифмической или
степенной функций. Например, неравенства для среднего
логарифмического \eqref{log1} имеют смысл уточнений теоремы Лагранжа о среднем значении для логарифмической функции. Действительно, по теореме о среднем получаем
$$
\ln(y)-\ln(x)=\frac{1}{a}(y-x), y>x,
$$
где для промежуточного значения $a$ из общей теоремы о среднем известна только оценка $x<a<y$.
Но очевидно равенство
$$
a=L(x,y)
$$
и соответствующая оценка, вытекающая из \eqref{log1}, даёт

$$
x<a=L(x,y)\leq A(x,y)= \frac{x+y}{2}.
$$
Таким образом, на самом деле среднее значение обязательно находится на левой половине исходного отрезка. А применение более точного неравенства Радо \eqref{log2} ещё сильнее прижимает возможный промежуток для среднего значения к левому концу исходного отрезка.

Большое число различных неравенств для среднего логарифмического и других средних, в том числе на комплексной плоскости,  можно найти в методичке для студентов \cite{Met}.

\subsection{Итерационные средние}

Рассмотрим итерационный процесс при
заданных стартовых значениях $x_0, y_0$  и паре
абстрактных средних $(M, N)$\ : \\
$$
x_{n+1} = M(x_n, y_n), \ y_{n+1} = N(x_n, y_n).
$$

В общем случае получаем некоторую динамическую систему на плоскости,
ко\-то\-рая имеет интересное асим\-пто\-ти\-чес\-кое по\-ве\-де\-ние
(динамику) при $n\rightarrow~\infty$, но очень сло\-жна для изучения
даже при простейшем выборе пары средних $(M, N)$.

\begin{definition} Пусть существует общий предел
последовательностей $x_n$ и $y_n$. Тогда он называется
\textit{итерационным средним} и обозначается
$$
\mu (M,N|\ x_0, y_0) = \mu(x_0, y_0) = \lim x_n = \lim y_n.
$$
\end{definition}

Самым известным итерационным средним является арифметико --
геометрическое, введенное Гауссом. Оно получается при выборе $M=M_1$
, $N=M_0$ \   и выражается по формуле
$$
\mu(M,N|\ x_0,y_0) = \frac{\frac{\pi}{2}\ x_0}{K\left(\sqrt{1-\left(
\frac{y_0}{x_0}\right)^{2}}\right)}, 0 < y_0 <x_0,
$$
где $K(x)$ - полный эллиптический интеграл Лежандра первого рода.
Теория итерационных средних в настоящее время активно развивается и
является, в частности, источником возникновения новых классов
специальных функций --- гипертрансцендентных функций. Непосредственно проверяется, что итерационное среднее
пары абстрактных средних также является новым абстрактным средним и
наследует соответствующие аксиомы 1) - 4).

На основе метода AGM (arithmetic--geometric mean)  вычисляются сами эллиптические интегралы, да и не только они. Сейчас это один из самых мощных методов  вычислительной математики. Заметную роль в пропаганде метода AGM и его развитии в конце 19 века  сыграл   Карл Борхард---один из учеников К.~Якоби. Именем Борхарда сейчас называется  класс итерационных средних. Более того, современные методы наиболее быстрого вычисления простейших элементарных функций, таких как экспонента, квадратный корень или синус с косинусом, также основаны на методе AGM. Это несколько напоминает ситуацию  с  современными оптимизированными алгоритмами для  обычного умножения чисел, основанными на быстром преобразовании Фурье (БПФ). В обоих случаях для вычисления элементарных выражений используются продвинутые математические теории и алгоритмы.

\subsection{Другие классы средних и неравенств для них}
 Кроме симметричных средних  используются и несимметричные. Самые
известные из них --- это  средние арифметическое и геометрическое с
неотрицательными весами $\alpha,\beta$ :
\begin{eqnarray}\label{ves1}
A_{\alpha,\beta}(x,y)=\alpha x + \beta y,\\
G_{\alpha,\beta}(x,y)=x^{\alpha} y^{\beta},~\alpha +\beta=1. \label{ves2}
\end{eqnarray}
Неравенство между ними---это знаменитое неравенство Янга (его часто называют неравенством Юнга). Отметим,
что его также можно уточнять методом средних, см. \cite{Sit1,Sit2}. Аналогично
вводится весовое среднее степенное произвольного порядка.

Несимметричные степенные средние также используются в Анализе,
например, в методах Рисса, Чезаро и Вороного обобщённого
суммирования рядов. Симметричные или несимметричные
весовые средние --- это также одно из основных понятий  теории
вероятностей и математической статистики.

\subsubsection {Средние Джини и Лемера}

Итальянский статистик Коррадо Джини ввёл названные его именем средние с двумя параметрами в 1938 г. по формулам
\begin{equation}\label{Jini}
Gi_{u,v}(x,y)=\left( \frac{x^u+y^u}{x^v+y^v}\right)^{\frac{1}{u-v}}, u\neq v,
\end{equation}
$$
Gi_{u,v}(x,y)=\exp \left(\frac{x^u\ln x+y^u \ln y}{x^u+y^u}\right), u=v\neq0,
$$
$$
Gi_{u,v}(x,y)=G(x,y), u=v=0.
$$

Основным предметом изучения Джини было имущественное неравенство, для его описания он и ввёл индекс Джини, основанный на предыдущих формулах. Это одно из основных понятий современной социальной статистики. Он также использовал и пропагандировал результаты своего соотечественника Вильфреда Парето. Парето первым провёл конкретные расчёты, из которых следовало, что в его время всего 20\% населения владело 80\% национального богатства, и сделал выводы, что организованное таким образом государство ждёт быстрая и неизбежная гибель.
Несмотря на то, что по запросу в интернете на множество Парето или индекс Джини будет выдано несколько десятков тысяч ссылок, найти информацию о них самих практически невозможно. Наверное, это связано с тем, что оба были теоретиками итальянского фашизма, в котором, видимо, видели способ возрождения недавно ставшей самостоятельной Италии. На идеях Парето воспитывался Муссолини, а Джини при нём занимал высокие посты и получал полную поддержку.

Важный частный случай средних Джини получается из \eqref{Jini} при $v=u-1$.
Эти средние были переоткрыты в  1971 г. Д.~Лемером и имеют вид
$$
Le_u(x,y)=\frac{x^{u+1}+y^{u+1}}{x^u+y^u}.
$$

Эти средние  в советское время чуть позже Лемера также переоткрыл школьный учитель Ю.~М.~Фирсов \cite{Kal}.
Разумеется, сейчас средние Джини и Лемера изучены вдоль и поперёк, в том числе для них получены теоремы типа Радо о сравнении сo средними других типов, см. \cite{BMV}.

\subsubsection{ Квазиарифметические средние}

\begin{definition} Пусть дана неотрицательная монотонная функция $f(x)$, набор неотрицательных чисел $x=(x_1,x_2,\cdots,x_n)$ и неотрицательных весов
$p=(p_1,p_2,\cdots,p_n)$. \textit{Квазиарифметическим средним} называется выражение
$$
K_p(x)=f^{-1}\left( \sum_{k=1}^n p_k f(x_k)\right), \sum_{k=1}^n p_k=1.
$$
\end{definition}
В частном случае $f(x)=x$ получаем обычное весовое среднее арифметическое, чем объясняется название. Это среднее стало широко известно после выхода монографии \cite{HLP}, в которой оно применялось для изучения обобщённой выпуклости относительно пары функций. Основные результаты для квазиарифметического среднего были получены в работах Колмогорова, Нагумо и Де Финетти. Они в существенном сводятся к доказательству двух основных свойств:\\

1. Квазиарифметические средние совпадают тогда и только тогда, когда порождающие их функции связаны линейным соотношением.\\

2. Непрерывное  квазиарифметическое среднее однородно тогда и только тогда, когда оно совпадает со средним степенным.\\

Данные средние встречаются во многих задачах и приложениях.  Они используются в некоторых работах  В.~П.~Маслова по применению нелинейных моделей в математической экономике \cite{Mas1,Mas2,Mas3}.

 Существуют и другие средние, задающиеся в явном виде:  Столярского (обобщённые средние Радо), Дирихле, итерационные средние Борхарда, интегральные
средние \- Карлсона в терминах гипергеометрических функций или в
виде симметрических интегралов и многие другие.
Полезным необычным средним является также понятие медианты двух
дробей, вводимое в теории рядов Фарея.

В заключение перечисления средних отметим, что некоторые
элементарные операции не выводят из их класса. Например, среднее от
пары средних - это вновь среднее. Несмотря на свою очевидность, это
формулировка очень мощного метода, порождающего в конкретных случаях
практически все известные процедуры построения новых средних из
простейших.

Таким образом, можно сделать основной для данного раздела вывод, что
\textit{существует значительное число конкретных примеров
абстрактных средних}. Они используются далее в качестве
"кирпичиков"{}, из которых собираются различные обобщения неравенств
Коши-Буняковского в интегральном и дискретном случаях.

Необходимо подчеркнуть фундаментальную роль средних степенных и
Радо, так как все остальные перечисленные выше известные средние
выражаются через них. Поэтому именно эти две шкалы средних играют
основную роль в данной работе.

Важность средних вытекает уже и из того очевидного факта, что самые
элементарные арифметические операции сложения и умножения выражаются
через простейшие средние:
$$
x+y=2 M_1 (x,y),\phantom{xxx} x\cdot y=( M_0(x,y))^2 ,
$$
так же, как и операции минимума и максимума.

Неравенства между средними --- это различные способы их
упорядочивания. Иногда используются достаточно экзотические варианты
такого упорядочивания, например, с использованием преобразования
Фурье и понятия положительно  определенных функций, которые тем не
менее находят полезные приложения, см. \cite{Fum}. Об использование средних при доказательстве неравенств для положительно  определённых функций см. \cite{PS}.

\subsection{ Неравенства для средних в комплексной плоскости }

Интересной темой является получение обобщений неравенств между средними в комплексной плоскости при аккуратном расставлении модулей  для комплексных величин. О простейших таких неравенствах см. \cite{Met}. Здесь только немного наметим эту тему.

Одна из основных структурных теорем алгебры запрещает сравнивать комплексные числа, так как утверждает существование единственного с точностью до изоморфизма упорядоченного числового поля, а это место уже занято действительными числами. Остаётся переходить к модулям.

Итак, пусть выполнено неравенство $f(x,y) \le g(x,y)$ для неотрицательных функций $f(x,y) \ge 0, g(x,y)\ge 0$ от \textit{действительных} переменных $x,y \ge 0$. Назовём  \textit{комплексификацией действительного неравенства} набор следующих неравенств  для \textit{комплексных} переменных $z,w\in \mathbb{C}$:
$$
|f(z,w)| \le |g(z,w)|,
$$
$$
|f(z,w)| \le g(|z|,|w|),
$$
$$
f(|z|,|w|) \le |g(z,w)|,
$$
$$
f(|z|,|w|) \le g(|z|,|w|).
$$
\begin{center}
\ldots \ldots \ldots
\end{center}

Рассматривать последнее из приведённых выше неравенств не имеет смысла, так как оно совпадает с исходным. На самом деле выписан далеко не полный комплект, что отражают многоточия, ведь модули можно вставить в трёх местах справа и трёх местах слева, итого получается $8 \times 8=64$ варианта, правда, некоторые из них отпадут. Например, неравенство $f(|z|,|w|) \le |g(|z|,w)|$ может быть добавлено в список, а вариант $f(z,|w|) \le |g(|z|,w)|$ следует забраковать, так как слева стоит комплексное число. Отдельно нужно рассматривать условия, когда неравенства переходят в равенства.

Приведённая схема описана в \cite{Met} на примере комплексификаций неравенства между средними арифметическим и геометрическим (автор раздавал подобные задачи своим студентам  для рисования графиков на первых появившихся компьютерах во Владивостокском политехе в конце 1980х -- начале 90х). Если сформулировать кратко, то справедлив неожиданный результат: на комплексной плоскости в подходящих переменных существует некоторая замкнутая кривая, напоминающая по форме улитку Паскаля, на которой в данном неравенстве достигается равенство, в её внешности (бесконечная область) по-прежнему среднее арифметическое больше геометрического, а внутри (конечная область) выполняется обратное неравенство.

Разумеется, приведённая схема подходит и для неравенств от одной переменной, и от нескольких, а не только от двух. Теперь можно применять введённую комплексификацию ко многим известным неравенствам. В принципе, процесс можно продолжить на две оставшиеся по теореме Гурвица нормированные гиперкомплексные алгебры---кватернионы и октавы.

Итак, дано исходное неравенство
$$
f(x,y)=\sqrt{xy} \le g(x,y)=(\frac{x+y}{2}), x\ge 0, y\ge 0.
$$
Корень нам не совсем подходит, поэтому, чтобы избежать многозначных функций, возведём это неравенство в квадрат. Его комплексификации выписываются тогда при $z,w\in \mathbb{C}$ так:
\begin{equation}\label{com4}
|zw| \le \left|\frac{z+w}{2}\right|^2,\qquad
|zw| \le {\left(\frac{|z|+|w|}{2}\right)}^2,\qquad
|zw| \le \left|\frac{|z|+w}{2}\right|^2.
\end{equation}
Для рассматриваемого неравенства весь полный боекомплект из 64 вариантов сводится к трём \eqref{com4}.

Мы хотим рассчитать и потом изобразить на графиках те множества, для которых выполняются неравенства \eqref{com4}. Сделать это в 4--мерном пространстве сложновато, поэтому упростим задачу и используем однородность в \eqref{com4}. Тогда для новой комплексной переменной $s=\frac{w}{z}$ получим
\begin{equation}\label{com5}
|s| \le \left|\frac{s+1}{2}\right|^2,\qquad
|s| \le {\left(\frac{|s|+1}{2}\right)}^2,\qquad
|s| \le \left|\frac{\frac{w}{|z|}+1}{2}\right|^2.
\end{equation}
Два первых неравенства теперь зависят от одной переменной, их мы и продолжим рассматривать. Что делать с третьим, которое опять зависит от двух комплексных переменных, непонятно, и мы его рассматривать не будем.

Исследуем первое из неравенств \eqref{com5}. После подстановки $s=x+iy$ и некоторых упрощений, получаем, что равенство в нём достигается на некоторой кривой 4 порядка с уравнением
\begin{equation}\label{com6}
x^4+y^4+2x^2y^2+4x^3+4xy^2-10x^2-14y^2+4x+1=0.
\end{equation}

C помощью Maple и некоторых расчётов можно сформулировать выводы.

\begin{theorem}
\begin{enumerate}
\item  На комплексной плоскости существует замкнутая ограниченная кривая, на которой в первом из неравенств \eqref{com5} достигается равенство, вне её выполнено неравенство, а в ограниченной области внутри кривой справедливо обратное неравенство. Указанная кривая определяется уравнением \eqref{com6}.\\

\item  Первое из неравенств \eqref{com5} для действительных $s$ справедливо не только для $s{>}0$, как следовало бы ожидать, а при $s\in[-\infty,-3-2\sqrt2]\cup[-3+2\sqrt2,\infty]$.\\

\item  Первое из неравенств \eqref{com5} для чисто мнимых $s$ выполняется при ${\rm Im}(s)\in [-\infty,-2-\sqrt3]\cup[-2+\sqrt3,2-\sqrt3]\cup[2+\sqrt3,\infty]$.\\

\item  Уравнения внешней и внутренней кривых в полярных координатах имеют вид $r=2-\cos{\phi} \pm \sqrt{(2-\cos{\phi})^2-1}$.
\end{enumerate}
\end{theorem}

\begin{figure}[h!]
	\begin{minipage}[h]{0.49\linewidth}
		\center{\includegraphics[scale=0.25]{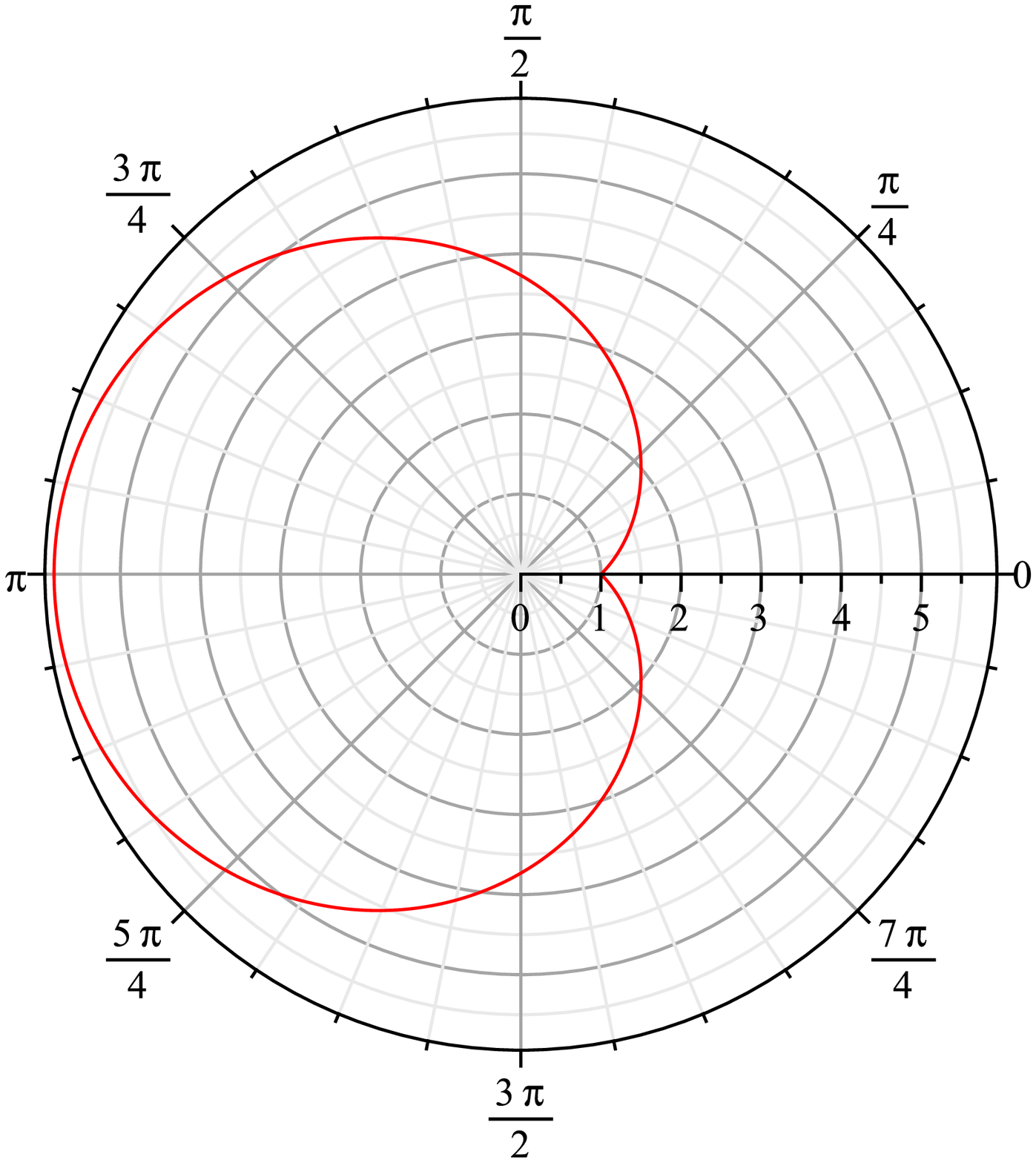} \\ a)}
	\end{minipage}
	\hfill
	\begin{minipage}[h!]{0.49\linewidth}
		\center{\includegraphics[scale=0.25]{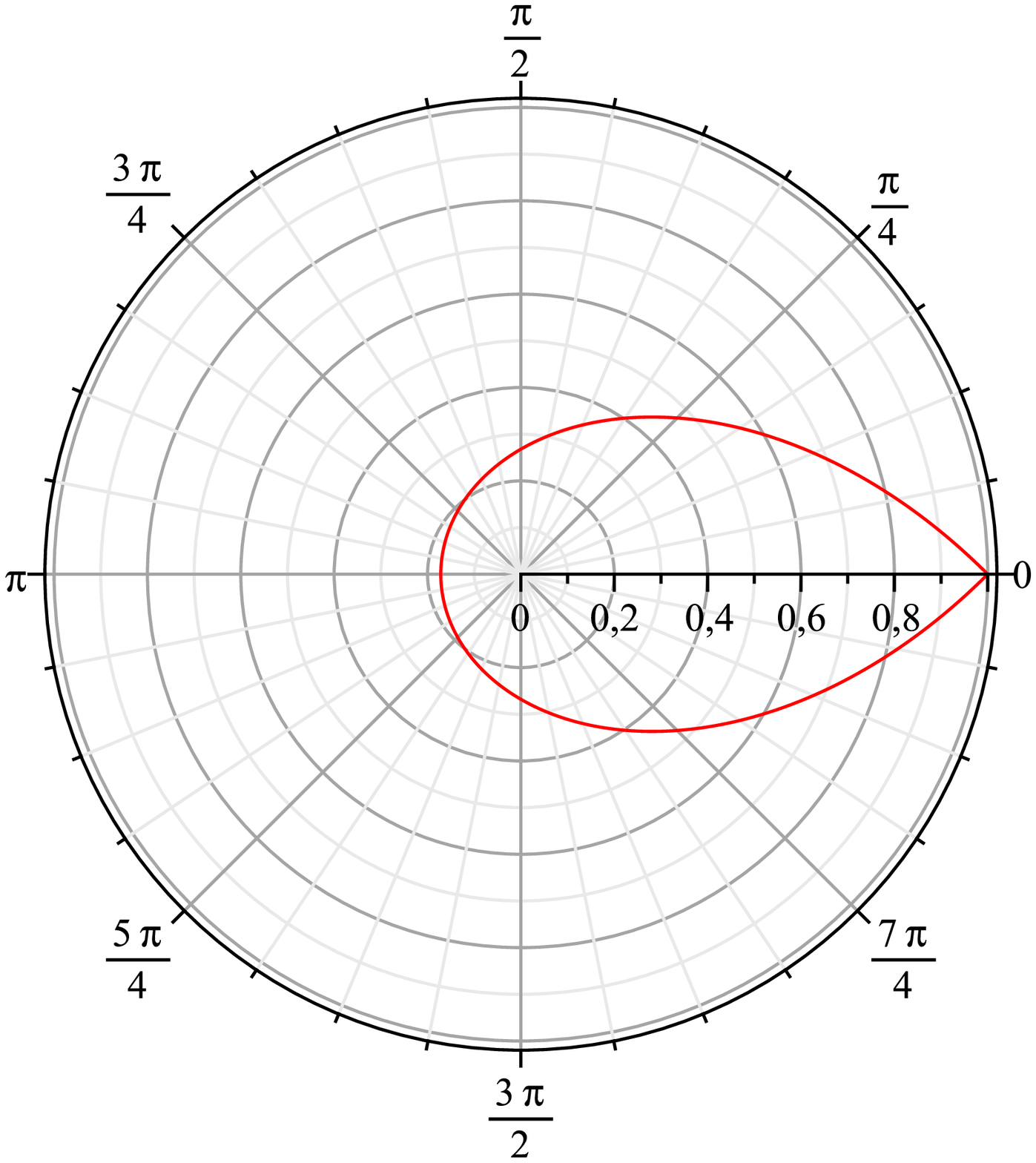} \\ b)}
	\end{minipage}
	\caption{Первая a) и вторая b) части кривой $r=2-\cos{\phi} \pm \sqrt{(2-\cos{\phi})^2-1}$ в полярных координатах.}
	\label{ris1}
\end{figure}

\begin{figure}[h!]
\begin{center}
\includegraphics[scale=0.3]{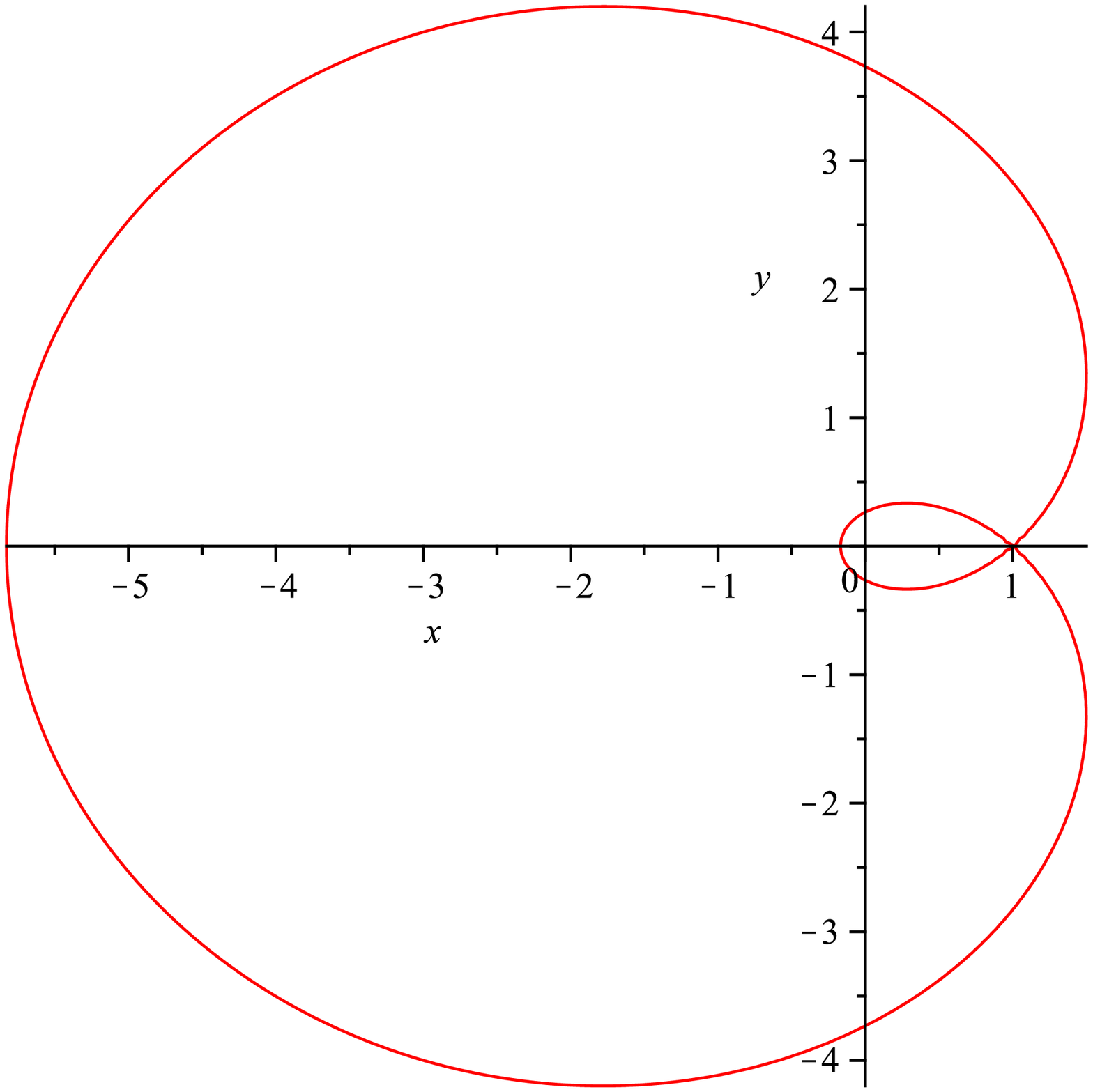}
\caption{Полный график  разделяющей кривой   в прямоугольных координатах.}\label{ris2}
\end{center}
\end{figure}

Приведём численные значения $-3-2\sqrt2\approx-5,8284, -3+2\sqrt2\approx-0,1716, -2-\sqrt3\approx-3,7321, -2+\sqrt3\approx-0,2679,  2-\sqrt3\approx0,2679, 2+\sqrt3\approx3,7321.$
Графики $r=2-\cos{\phi} \pm \sqrt{(2-\cos{\phi})^2-1}$ сначала отдельно двух кусков в полярных координатах, а затем всей кривой в прямоугольных координатах приведены  на Рис. \ref{ris1} и \ref{ris2}.
Получившиеся кривые похожи на некоторые известные кривые 4 порядка,
а именно, улитки Паскаля, кардиоиды и овалы Декарта, но к ним не сводятся.

Что же касается второго из неравенств \eqref{com5}, то с ним всё ясно. Это неравенство только для одного радиуса и оно выполняется на границе и вне единичного круга и не выполняется внутри его.

Таким образом, указанная процедура комплексификации неравенств может рассматриваться как источник новых плоских кривых. Можно назвать их разделяющими для комплексификаций соответствующего неравенства.

Следует сказать, что неравенства для средних в комплексной плоскости --- это не только интересная забава.  В знаменитой работе В.~Бекнера \cite{Bec}, в которой получены точные значения норм преобразования Фурье в
шкале пространств $L_p$, используются две основные конструкции:
дробное (или квадратичное) преобразование Фурье и некоторое неравенство для средних степенных  в комплексной
плоскости. Точную константу, которую получил Бекнер, искали и многие до него. Наибольших успехов добился наш математик К.~И.~Бабенко, решивший задачу в частном случае $p=2^n$, поэтому оправдано название "неравенство Бабенко--Бекнера". Ему принадлежит идея применять в доказательстве специальное квадратичное преобразование Фурье. Вторая идея принадлежит Бекнеру, он свёл оставшуюся часть задачи к некоторому неравенству между степенными средними в комплексной плоскости. Нужно было доказать, что это неравенство выполняется на определённом отрезке мнимой оси, что Бекнеру и удалось сделать.

\subsection{Характеризация произвольных средних для двух чисел.\\ Связь средних и энтропии}

~

\medskip
Можно пойти дальше и получить  описание всех
"хороших"\ средних от двух переменных.

\begin{theorem} \label{th}\textit{Произвольное абстрактное среднее},
обладающее свойствами несмещённости, однородности, симметричности,
монотонности и непрерывности по обоим аргументам
\textit{представляется в виде}
\begin{equation}\label{h}
M (x,y)= (x+y)h(\ln \frac{y}{x}),
\end{equation}
 где  $h(t)$ является определенной на всей оси непрерывной
 четной функцией, удовлетворяющей следующим условиям при $t\geq
0$:
\begin{eqnarray}
h(0)=\frac{1}{2} ,\\
\frac{e^{t_1}(e^{t_2}+1)}{e^{t_2}(e^{t_1}+1)}\ \leq \
\frac{h(t_1)}{h(t_2)}& \leq & \frac{(e^{t_2}+1)}{(e^{t_1}+1)}, \
t_1\leq t_2.
\end{eqnarray}
\textit{Справедливо и обратное:} каждой функции $h(t)$ с указанным
набором свойств соответствует по формуле \eqref{h} несмещённое,
однородное, симметричное, монотонное и непрерывное по обоим
аргументам среднее.
\end{theorem}

Введенная функция $h(t)$ связана с функцией однородности
соответствующего среднего $M (x,y)$. Отметим, что свойство чётности для функции $h(x)$ равносильно симметричности среднего $M (x,y)$ в \eqref{h}.

Полученное полное описание средних значений позволяет дать
исчерпывающее описание функций энтропии от двух переменных, что
является важным для теории информации и других приложений.
Действительно, определим обобщённую функцию энтропии $ H$  через
произвольное среднее по  формуле
 \begin{equation} \label{entr}
  H=-\ln( M( p_1,p_2 )).
  \end{equation}
  Это определение не надумано, а отражает известные конкретные функции энтропии, которые общеприняты в термодинамике и информатике. Получается,  что при таком подходе с использованием теоремы \ref{th}
нами получено полное внутреннее описание всех двумерных функций
энтропии. В приведённых обозначениях $p_1,p_2$ --- это вероятности значений некоторой дискретной случайной величины с функцией распределения (средним)  $M( p_1,p_2 )$, подчиняющиеся соотношению $p_1+p_2=1$. При этом выбору специального весового среднего геометрического \eqref{ves2} соответствует энтропия по Шеннону:
$$
M(p_1,p_2)=G_{p_1,p_2}(p_1,p_2)=p_1^{p_1}\cdot p_2^{p_2}, H_{Sh}(p_1,p_2)=-p_1\ln(p_1) - p_2\ln(p_2),\  p_1+p_2=1;
$$
а  выбору  среднего Джини \eqref{Jini} соответствует энтропия по Реньи:
$$
M(p_1,p_2)=Gi_\alpha(p_1,p_2)=\left(\frac{1}{p_1^\alpha+p_2^\alpha}\right)^{\frac{1}{1-\alpha}}=
\left(\frac{p_1+p_2}{p_1^\alpha+p_2^\alpha}\right)^{\frac{1}{1-\alpha}}, p_1+p_2=1,
$$
$$
H_{Re,\alpha}(p_1,p_2)=\frac{1}{1-\alpha}\ln(p_1^\alpha+p_2^\alpha), \ p_1+p_2=1.
$$
Эти соотношения можно продолжить на любое число вероятностей, а также на другие варианты средних и соответствующих энтропий.

Приведём следующее из предыдущей теоремы \ref{th}  полное описание всех функций энтропии от двух переменных.

\begin{theorem} Произвольная функция энтропии, определяемая соотношением  \eqref{entr},   может быть получена по формуле
$$
H(p_1,p_2)=-\ln(p_1+p_2)-\ln \left(h(\ln \frac{p_2}{p_1})\right),
$$
с некоторой функцией  $h(t)$, где  $h(t)$ является определенной на всей оси непрерывной  четной функцией, удовлетворяющей следующим условиям при $t\ge
0$:
$$
h(0)=\frac{1}{2} ,\\
$$
$$
\frac{e^{t_1}(e^{t_2}+1)}{e^{t_2}(e^{t_1}+1)}\ \le \
\frac{h(t_1)}{h(t_2)} \le  \frac{(e^{t_2}+1)}{(e^{t_1}+1)}, \
t_1\le t_2.
$$
\end{theorem}

По-видимому, предложенный способ описывает только функции энтропии, удовлетворяющие условию аддитивности. Поэтому он охватывает энтропии по Шеннону и Реньи, но не описывает ряд других известных функций энтропии, таких, например, как энтропии Цаллиса, Чисара и некоторые другие.

\subsection{ Все классические неравенства --- это неравенства для средних!}

~

\medskip
Действительно, мы уже обсуждали, что неравенство Янга---это неравенство между весовыми геометрическим и арифметическим средними. Непоследственно в терминах средних формулируется известная теорема Ляпунова, см. \cite{HLP}.

Далее, рассмотрим пару непрерывных (только чтобы не усложнять ситуацию) на отрезке $[a,b]$ неотрицательных функций $f(x),g(x)$. Введём ещё функционалы
\begin{equation} \label{n1}
J(f)=\int_a^bf(x)\ dx, K(f)=\sqrt{J(f)}=\sqrt{\int_a^bf(x)\ dx}.
\end{equation}
Напомним обозначения для средних арифметического, геометрического и весового геометрического: $A,G,G_{u,v}$.

Тогда интегральное неравенство Коши--Буняковского во введённых обозначениях записывается так:
\begin{equation} \label{n2}
J(G^2(f,g))\le G(J(f^2),J(g^2)).
\end{equation}

Неравенство Гёльдера (точнее Роджерса-Гёльдера-Рисса) во введённых обозначениях записывается так:
\begin{equation} \label{n3}
J(G^2(f,g))\le G_{\frac{1}{p},\frac{1}{q}}(J(f^p),J(g^q)),\frac{1}{p}+\frac{1}{q}=1 .
\end{equation}

Неравенство Минковского во введённых обозначениях записывается так:
\begin{equation} \label{n4}
K(A^2(f,g))\le A(K(f^2),K(g^2)) .
\end{equation}

В терминах средних формулируются и оставшиеся теоремы из классического набора неравенств: Чебышёва, Иенсена и Артина. Теоретически доказано, что все перечисленные в этом пункте неравенства эквивалентны, то есть известны способы доказательства, что они следуют друг из друга \cite{MO2}. (Тут важно именно то, что известны явные доказательства, с помощью которых неравенства выводятся друг из друга, а с общей точки зрения, как известно из курса математической логики, все истинные утверждения эквивалентны!).  Но некоторые из указанных следствий конкретно и явно проследить крайне сложно, несмотря на то, что теоретически доказано, что они существуют.

Ничего нового в утверждениях \eqref{n1}--\eqref{n4} по существу не содержится. Но они приводят к такому очевидному вопросу, как найти необходимые и достаточные условия для выполнения каждого из неравенств \eqref{n1}--\eqref{n4}, заменив конкретные средние на произвольные абстрактные средние, а функционалы $J,K$ на пару произвольных функционалов. Представляется, что это очень трудные задачи. Каждую можно попробовать начать решать по частям, заменяя на произвольную величину что--то одно.

Далее, рассмотрим некоторое нормированное пространство, на элементах которого определены все встречающиеся ниже средние. Тогда неравенство треугольника можно записать с использованием среднего арифметического так:
\begin{equation}\label{tr1}
 \|A(x,y)\| \le A(\|x\|,\|y\|).
\end{equation}
Можно поставить задачу об описании всех допустимых средних $A$, для которых выполнено неравенство треугольника \eqref{tr1} в данном нормированном пространстве.

Более общая задача состоит в описании всех пар допустимых средних $M,N$, для которых выполнено обобщённое неравенство треугольника
\begin{equation}\label{tr2}
\|M(x,y)\| \le N(\|x\|,\|y\|)
\end{equation}
для некоторой нормы $\|\cdot\|$, а также описать все нормы, для которых неравенство \eqref{tr2} возможно при некоторых парах средних $M,N$ с явным указанием таких пар.

\subsection{Два интересных приложения неравенств для средних: результаты Дю--Плесси и неравенства для тета--функций Якоби.}

Эти приложения малоизвестны, но великолепно отражают силу неравенств для средних. Первое взято из знаменитой монографии \cite{SKM} по дробному интегродифференцированию (Красная Книга!), а второе --- из  малоизвестного обзора  по неравенствам для тета--функций Якоби \cite{MS}.

\subsubsection{Оценки Дю--Плесси} Неравенства С.~Л.~Соболева в простейшей форме---это оценки сингулярного оператора
$$
\int_{R^n} \frac{f(y)\ dy}{|x-y|^{n-\alpha}}
$$
при подходящих $\alpha$ в нормах пространств Лебега  $L_p(R^n)$ (одного или двух), в знаменателе стоит расстояние между двумя векторами в степени. В одномерном случае есть знаменитая теорема Харди--Литтвульда, обеспечивающая нужную оценку. Затем появились десятки работ, в которых разными способами неравенство Соболева переносилось на многомерный случай.

Далее просто процитируем энциклопедическую монографию по дробному исчислению \cite{SKM}. "Приём N.~du~Plessis... сведения теоремы Соболева к случаю $n=1$ состоит в наблюдении, что \ldots $|x|^n \ge n^{\frac{n}{2}}\prod_{k=1}^n |x_k|$, что вытекает из того факта, что среднее арифметическое не меньше, чем среднее геометрическое. Поэтому
\begin{equation}\label{Sob}
\int_{R^n} \frac{|f(y)|\ dy}{|x-y|^{n-\alpha}}\le C
\int_{R^n} \frac{|f(y)|\ dy}
{\prod_{k=1}^n |x_k-y_k|^
{1-
\frac{\alpha}{n}}},
\end{equation}
и при оценке достаточно располагать одномерной теоремой Харди--Литтвульда, применив её по каждой переменной". Таким образом, благодаря Дю Плесси, догадавшемуся к данному интегралу применить неравенство между средним квадратичным (или арифметическим) и средним геометрическим, стала ненужной огромная куча работ, разросшихся по данной тематике. Но стоит отметить, что при таком подходе невозможно установить точные постоянные в неравенстве \eqref{Sob}, а также рассматривать неравенства в областях общего вида..

\subsubsection{ Неравенства для тета--функций Якоби} Рассмотрим третью тета--функцию Якоби, которая определяется при помощи быстро сходящегося  тригонометрического ряда по формуле
\begin{equation}\label{teta}
\vartheta_3(z,q)=1+2 \sum_{k=1}^\infty q^{k^2} \cos(2kz),
\end{equation}
где $0<q<1, 0 \le z \le \pi$. Эта функция была введена К.~Якоби в 1829 г..
Данная тета--функция является целой по $z$ и периодической с периодом $\pi$ на действительной оси. Нетрудно показать, что функция имеет один положительный минимум посередине периода, при этом величина этого минимума $m(q)$ очень мала, если $q\to 1$. Например, $m(0,999)\approx 10^{-1069}$. Рассмотрим задачу об оценке этого минимума сверху, так как снизу хорошая оценка---это просто ноль.

\begin{theorem} Для величины минимума тета--функции Якоби справедлива оценка
\begin{equation}\label{min1}
  m(q)=\min \vartheta_3(z,q)=\vartheta_3(\frac{\pi}{2},q)\le
\exp
\left(
-\frac{q^2+2q}{1-q^2}
\right).
\end{equation}
\end{theorem}

Заметим, что при граничных значениях $q=0$ и $q=1$ в применённом неравенстве о средних достигается знак равенства. Следовательно, наши оценки точны вблизи этих значений  $q=0$ и, что особенно важно, при $q=1$. Эти значения наиболее интересны, так как функция имеет в пределе особенность и не имеет там приличной асимптотики на комплексной плоскости вблизи особой точки.

Расчёт показывает, что при $q=0,999$ неравенство \eqref{min1} даёт $m(0,999)\le 1,2*10^{-651}$. Неплохо, мы  получили около 60{\%} от фактического порядка минимума, как показывает компьютерный расчёт. Дальнейшие уточнения полученной оценки и анализ её точности см. в \cite{MS}. Понятно, что подобные вычисления в компьютерном пакете MATHEMATICA необходимо вести с использованием чисел с огромным числом разрядов.

Как же всё--таки можно попробовать получать более точные оценки, которые приближались бы для рассматриваемого нами минимума к истинному порядку этой величины?  Один из путей, который представляется возможным, это использовать различные уточнения неравенства между арифметическим и геометрическим средними. Таких уточнений известно в литературе по неравенствам немало, возможно какие-то из них и приведут к намеченной цели, но это пока дело будущего. В том числе для этой цели можно использовать некоторые результаты этой работы.

\section{ Обобщения дискретных неравенств Коши--Буняковского}

Итак, рассмотрим дискретное неравенство Коши--Буняковского
\begin{equation}\label{D1}
\left(\sum_{k=1}^nx_k\cdot y_k\right)^2 \le \left(\sum_{k=1}^nx_k^2\right)\cdot\left(\sum_{k=1}^ny_k^2\right).
\end{equation}

В \cite{HLP} приведено его интересное обобщение  в форме
\begin{equation}\label{D2}
\left(\sum_{k=1}^nx_k\cdot y_k\right)^2 \le
\left(\sum_{k=1}^n(x_k^2+y_k^2)\right)\cdot\left(
\sum_{k=1}^n\frac{x_k^2 y_k^2}{x_k^2+y_k^2}\right)
\le
\end{equation}
$$
\le\left(\sum_{k=1}^nx_k^2\right)\cdot\left(\sum_{k=1}^ny_k^2\right).
$$

Это  известное  неравенство,  которое было доказано Милном
в 1925 году при исследовании задачи из астрономии
по вычислению коэффициента звездного поглощения, и с тех пор
неоднократно переоткрывалось. Неравенство
Милна явилось для меня одной из отправных точек для получения изложенных далее результатов по обобщению неравенства Коши--Буняковского.

Другим подобным результатом является неравенство, доказанное Колбо в  1965 г., см. \cite{Dr1}.
$$
\left(\sum_{k=1}^nx_k\cdot y_k\right)^2 \le
\sum_{k=1}^nx_k^{1+\alpha} y_k^{1-\alpha}
\sum_{k=1}^nx_k^{1-\alpha} y_k^{1+\alpha}\le
\left(\sum_{k=1}^nx_k^2\right)\cdot\left(\sum_{k=1}^ny_k^2\right),
$$
справедливое при $\alpha\in[0,1]$.

Достаточно очевидно, что приведённые выше обобщения неравенства Коши--Буняковского, принадлежащие Милну и Колбо, имеют общую природу. Она проявляется в следующем результате, который называется теоремой Дейкина--Элиезера--Карлица ( теорема CDE).

\begin{theorem}
Для того, чтобы для любых двух последовательностей неотрицательных чисел было выполнено неравенство
\begin{equation} \label{CDE1}
\left(\sum_{k=1}^nx_k\cdot y_k\right)^2 \le
\sum_{k=1}^n f(x_k,y_k)\sum_{k=1}^n g(x_k,y_k)\le
\end{equation}
\begin{equation}\label{CDE11}
\left(\sum_{k=1}^nx_k^2\right)\cdot\left(\sum_{k=1}^ny_k^2\right),
\end{equation}
необходимо и достаточно, чтобы пара функций $f(x,y),g(x,y)$ удовлетворяла следующим трём условиям
\begin{equation}\label{CDE2}
1) \  f(x,y) g(x,y)=x^2 y^2,
\end{equation}
\begin{equation}
2) \ f(\lambda x,\lambda y)=\lambda^2 f(x,y),
\end{equation}
\begin{equation}\label{CDE3}
3) \ \frac{y f(x,1)}{x f(y,1)}+\frac{x f(y,1)}{y f(x,1)}\le \frac{x}{y}+\frac{y}{x}.
\end{equation}
\end{theorem}

Неравенство Милна является частным случаем теоремы CDE при выборе $f(x,y)=x^2+y^2$, а неравенство Колбо при выборе $f(x,y)=x^{1+\alpha} y^{1-\alpha}, \alpha\in[0,1]$.

Если проанализировать доказательство теоремы CDE, то становится ясным, что в достаточной части авторы переоткрывают аналог тождества Лагранжа, который существует для неравенства Чебышёва. Поэтому доказательство достаточности лучше сразу провести с применением неравенства Чебышёва, что и проще и не затемняет сути дела. Замечательным является найденный авторами способ доказательства необходимости теоремы и обнаруженное ими условие \eqref{CDE3}. Это гибридное условие на функцию, записываемое единой простой формулой, из которого следует одновременно её однородность и монотонность, аналогично тому как из известного условия Годуновой--Левина следует одновременно монотонность и выпуклость, см., например, \cite{HLP}. На мой взгляд роль таких гибридных условий, объединяющих в одной формуле несколько стандартных, пока в Анализе недооценена, их значение ещё предстоит раскрыть.

Теперь сделаем следующий важный шаг. Перепишем неравенства Милна и Колбо в терминах средних значений. Начнём с первого из них.
\begin{equation}\label{D3}
\left(\sum_{k=1}^nx_k\cdot y_k\right)^2 \le
\left(\sum_{k=1}^n {(M_2(x_k,y_k))}^2\right)\cdot\
\left(\sum_{k=1}^n {({M_2^*}(x_k,y_k))}^2\right)
\le
\end{equation}
\begin{equation}\label{D3}
\le\left(\sum_{k=1}^nx_k^2\right)\cdot\left(\sum_{k=1}^ny_k^2\right),
\end{equation}
где $M_2=Q$---среднее квадратичное, а $M_2^*=M_{-2}$---его дополнительное среднее, см. определение  \eqref{1}. Это предлагаемая мной запись неравенства Милна через среднее квадратичное.

Далее, запишем
\begin{equation}
\left(\sum_{k=1}^nx_k\cdot y_k\right)^2 \le
\left(\sum_{k=1}^n {(G_{\frac{1+\alpha}{2},\frac{1-\alpha}{2}}(x_k,y_k))}^2\right)\cdot\
\end{equation}
\begin{equation}\label{D4}
\cdot\left(\sum_{k=1}^n {({G^*}_{\frac{1+\alpha}{2},\frac{1-\alpha}{2}}(x_k,y_k))}^2\right)
\le\left(\sum_{k=1}^nx_k^2\right)\cdot\left(\sum_{k=1}^ny_k^2\right),
\end{equation}
где $\alpha\in[0,1]$, $G_{\frac{1+\alpha}{2},\frac{1-\alpha}{2}}$---весовое среднее геометрическое, ${G^*}_{\frac{1+\alpha}{2},\frac{1-\alpha}{2}}=
G_{\frac{1-\alpha}{2},\frac{1+\alpha}{2}}$---его дополнительное среднее. Это предлагаемая мной запись неравенства Колбо через весовое среднее геометрическое.

Теперь мы готовы записать теорему Дэйкина--Элиезера--Карлица полностью в терминах средних, начав соединение основных тем данного обзора.

\begin{theorem} (теорема CDE в терминах средних значений)
Пусть дано абстрактное несмещённое, однородное и монотонное по обоим аргументам среднее $M(x,y)$, $M^*(x,y)$---его дополнительное  среднее, см. определение  \eqref{1}. Тогда справедливо следующее обобщение дискретного неравенства Коши--Буняковского
\begin{equation}\label{D5}
\left(\sum_{k=1}^nx_k\cdot y_k\right)^2 \le
\left(\sum_{k=1}^n {(M(x_k,y_k))}^2\right)\cdot\
\left(\sum_{k=1}^n {({M^*}(x_k,y_k))}^2\right)
\le
\end{equation}
\begin{equation}\label{D55}
\le\left(\sum_{k=1}^nx_k^2\right)\cdot\left(\sum_{k=1}^ny_k^2\right).
\end{equation}

Верно и обратное: любое обобщение неравенства Коши--Буняковского из теоремы
Дэйкина--Элиезера--Карлица (CDE) в форме \eqref{CDE1}--\eqref{CDE11} может быть записано в виде \eqref{D5}--\eqref{D55} с некоторым средним $M$,  удовлетворяющим перечисленным выше свойствам (несмещённое, однородное и монотонное по обоим аргументам), и его дополнительным средним $M^*$.
\end{theorem}

Отметим, что условие симметричности среднего в этой теореме не требуется, что допускает такие случаи, как неравенство Колбо.

По моему мнению, переформулировка теоремы
CDE в терминах средних делает ее гораздо более понятной и
прозрачной, а также снабжает эту теорему больщим числом конкретных
примеров с использованием теории  средних величин в дополнение к скромному набору
из двух примеров (неравенства Милна и Колбо) во всех  известных монографиях, см., например, \cite{Dr1,Ste}.
Становится понятным, зачем мы описывали большое множество средних в предыдущей части  данного обзора --- с их помощью теперь можно повыписывать также большое число обобщений неравенства Коши--Буняковского совершенно различного вида. Важной задачей является также сравнение этих обобщений между собой, когда это возможно. Несмотря на соблазн это всё проделать, я здесь остановлюсь, и мы перейдём к рассмотрению интегрального случая, который представляется более интересным. К приложениям приведённых результатов для дискретного случая мы  вернёмся.

\section{Обобщения интегральных неравенств Коши--Буняковского}

\subsection{Интегральный аналог теоремы Дейкина--Элиезера--Карлица в терминах средних}

Перед формулировкой дальнейших результатов примем простейшие соглашение о том, что все далее встречающиеся функции являются непрерывными, зависят од одной  переменной, промежутки интегрирования являются конечными, а интегралы понимаются в смысле Римана и существуют. Эти условия далее в формулировках будем опускать, они предполагаются выполненными. От всех перечисленных ограничений можно избавиться, но для этого нужны отдельные рассмотрения. Мы также не будем исследовать важный вопрос, когда в рассматриваемых неравенствах достигается  равенство.

Одним из основных результатов теории, развитой  автором в \cite{Sit1,Sit3,Sit4,Sit2}, является следующее уточнение интегрального неравенства Коши--Буняковского.

\begin{theorem} \label{T1}
Пусть  $M$ - произвольное абстрактное среднее, удовлетворяющее перечисленным выше в определении 1 аксиомам 1--4, $M^*$ - дополнительное к нему среднее. Тогда справедливо обобщение
интегрального неравенства Коши--Буняковского вида
\begin{equation}\label{6}
  \left(\int^b_a f(x)g(x)\,dx\right)^2 \le
\int^b_a (M(f,g))^2\,dx\cdot
 \int^b_a(M^*(f,g))^2\,dx \le
\end{equation}
\begin{equation}\label{7}
 \le\int^b_a(f(x))^2\,dx \cdot
\int^b_a(g(x))^2\,dx , \phantom{aaaa}
\end{equation}
\end{theorem}

Таким образом, мы получаем переформулированный согласно работам \cite{Sit1,Sit3,Sit4,Sit2} в терминах средних результат, что справедлив интегральный аналог \textit{достаточной}
части  теоремы Карлица--Дэйкина--Элиезера (CDE, см. \cite{Sit1,Sit3,Sit4,Sit2}) для дискретного случая. В работах \cite{Sit1,Sit3,Sit4,Sit2} также доказан неожиданный результат, что в интегральном случае аналог \textit{необходимой} части этой известной теоремы не выполняется.
Анализ доказательства теоремы CDE выявляет то место, в котором
рассуждения для случая сумм нельзя повторить для случая интегралов.
В теореме CDE используется возможность составить сумму из двух
слагаемых, при этом получается замечательное условие в виде одной
формулы, которое является \textit{одновременным обобщением условий
монотонности и однородности}, см. об этом условии \eqref{CDE3}
 выше. Но для интегралов нет условия,
аналогичному взятию суммы из двух первых слагаемых.

Сделаем важное для дальнейшего замечание, что первое из неравенств  \eqref{6} очевидно,
так как оно само сводится к обычному неравенству Коши -- Буняковского
ввиду определения дополнительного среднего и вытекающего из него равенства
$$
M(f,g) M^*(f,g)=fg.
$$
Таким образом,
содержательным в теореме \ref{T1} является правое неравенство между
произведениями интегралов. Однако, и левое очевидное неравенство из теоремы
\ref{T1} может быть с успехом использовано (см. далее об итерационных процедурах уточнения оценок для специальных функций).

Теперь выбор любых различных средних и дополнительных к ним при помощи теоремы \ref{T1} порождает разнообразные обобщения вида (\ref{6}) интегрального неравенства Коши--Буняковского.
В частности, можно выписать варианты обобщений со степенными средними и средними Радо.
При выборе арифметико--геометрического среднего Гаусса получаем

\begin{corollary} Справедливо следующее уточнение неравенства
Коши -- Буняковского:
 \begin{eqnarray*}
 \left(
\int^b_af(x)g(x)dx\right)^2\le\int^b_a{\left[\frac{\max(f,g)}
{K\left(\sqrt{1-\left(\frac{\min(f,g)}
{\max(f,g)}\right)^2}\right)}\right]}^2dx \cdot\\
\cdot\int^b_a ( \min(f,g) )^2 \Biggl(
K\left(\sqrt{1-\left(\frac{\min(f,g)} {\max(f,g)}\right)^2}\right)
\Biggr) ^2\,dx \le\int^b_af^2\,dx\int^b_ag^2\,dx,
\end{eqnarray*}
где $K(x)$ есть полный эллиптический интеграл Лежандра 1 рода.
\end{corollary}

Отметим совершенно экзотический характер последнего неравенства: это
неравенство между \textit{произвольными} функциями, но которые стоят
под знаком \textit{конкретной специальной} функции --- эллиптического
интеграла Лежандра!

\subsection{Обобщение интегрального неравенства Коши--Буняковского в терминах максимума и минимума }

Рассмотрим теперь случай, когда среднее сводится к максимуму, а его дополнительное среднее---к минимуму двух функций (или наоборот). В результате получается

\begin{theorem}\label{mm} Пусть функции $f(x,y), g(x,y)$ являются неотрицательными на отрезке $[a,b]$. Тогда справедливо следующее уточнение
неравенства Коши--Буняковского
\begin{equation}\label{10}
\left(\int^b_af(x)g(x)dx\right)^2\le\int^b_a[\max(f,g)]^2dx\cdot
\int^b_a[\min(f,g)]^2dx\le
\end{equation}
$$
\le\int^b_af^2(x)dx\cdot \int^b_ag^2(x)dx.
\phantom{aaaaaaaaaaa}
$$
\end{theorem}

При первом взгляде на  неравенство (\ref{10}) кажется, что средняя часть
должна сводится к одной из крайних. Однако, это не так. То, что все
три части в (\ref{10}) могут быть различны, доказывает пример "конвертика" \ , образованного графиками функций \  $f(x)=x$, $g(x)=1-x$ при выборе пределов интегрирования $a=0$, $b=1$.
Вычисления показывают, что в этом случае (\ref{10}) сводится к строгим
неравенствам $\frac{1}{36}< \frac{7}{144}< \frac{1}{9}$. \smallskip

\begin{figure}[h!]
\center{\includegraphics[scale=0.4]{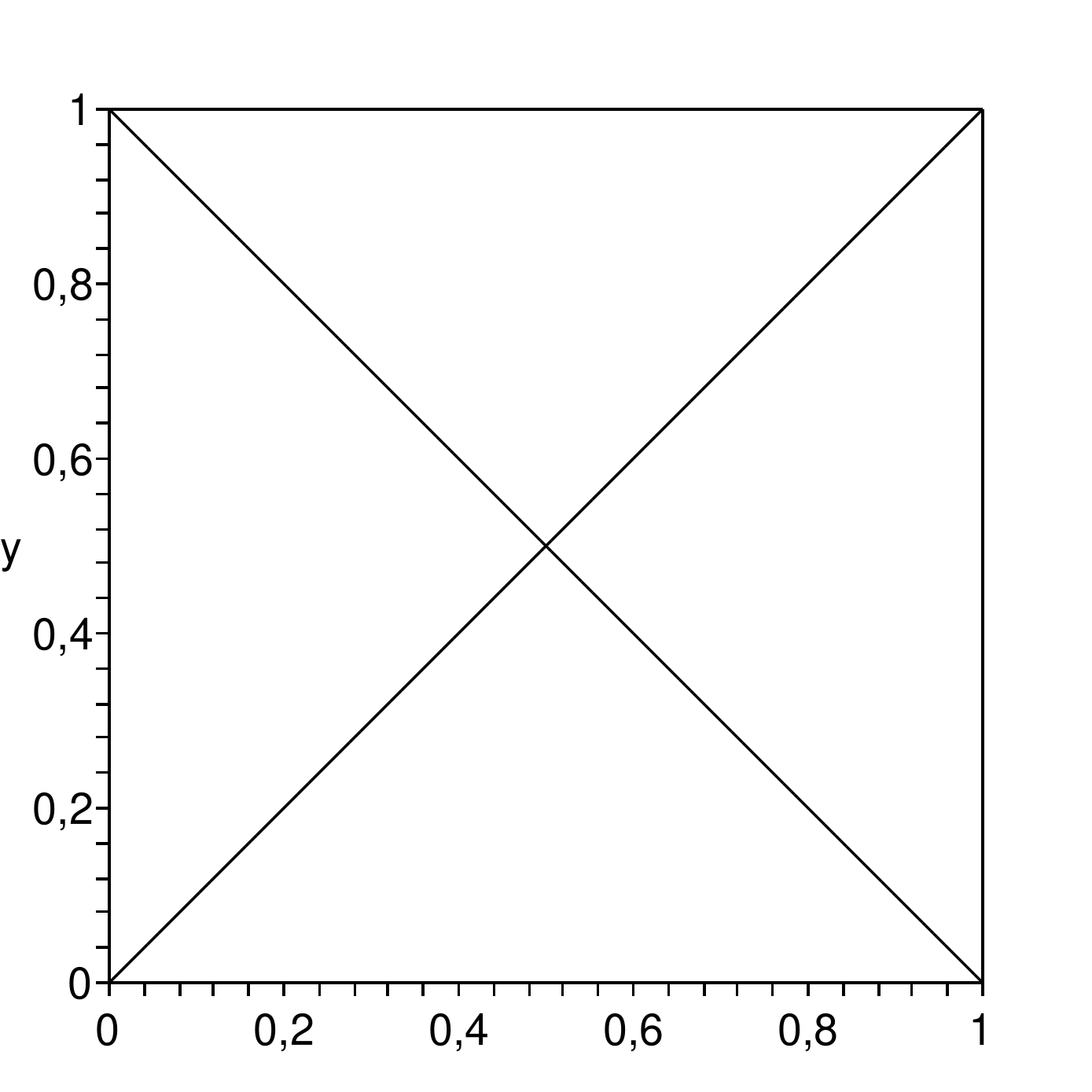}}
\caption{Пример "конвертика": $f(x)=x$, $g(x)=1-x$, $a=0$, $b=1$.}
\end{figure}

Последнее неравенство подробно изучено в работе \cite{Sit4}. В том числе там показано, что оно выполняется и при отказе от положительности функций, достаточно выполнения условия $f(x)+g(x) \geq 0$, что неожиданно. Приведём некоторые из этих результатов подробнее.

Начнём с доказательства теоремы \ref{mm}. Отметим, что она является частным случаем приведённой выше теоремы \ref{T1}. Первое прямое  доказательство было получено Анатолием Кореновским. Мы приводим доказательство автора, которое использует явные формулы  для максимума и минимума.

\textbf{ Доказательство.}

Известны следующие формулы  для максимума и минимума, которые являются следствиями теоремы Виета для квадратного трёхчлена:
$$
m=\frac{a+b-|a-b|}{2}, M=\frac{a+b+|a-b|}{2},\ m=\min(a,b), M=\max(a,b).
$$

Чтобы перейти к доказательству неравенства \eqref{10},  обозначим
$$
I(f)=\int_a^b f(x)\,dx.
$$
Тогда левая часть неравенства \eqref{10}  равна
$$
I((\min(f,g))^2)I((\max(f,g))^2)=
$$
$$
=I\left({\left(\frac{f+g-|f-g|}{2}\right)^2}\right)\cdot
I\left({\left(\frac{f+g+|f-g|}{2}\right)^2}\right)=
$$
$$
=\frac{1}{16} I\left((f+g)^2-2 (f+g) |f-g|+(f-g)^2\right)\cdot
$$
$$
 I\left((f+g)^2+2 (f+g) |f-g|+(f-g)^2\right)=
 $$
 $$
=\frac{1}{16} I\left(2f^2+2g^2-2(f+g) |f-g|\right)\cdot
 I\left(2f^2+2g^2+2 (f+g) |f-g|\right)=
 $$
 $$
=\frac{1}{4} I\left(f^2+g^2-(f+g)|f-g|\right)\cdot
 I\left(f^2+g^2+2(f+g)|f-g|\right).
$$
После очевидных преобразований отсюда следует
$$
\frac{1}{4}\left[(I(f^2))^2+(I(g^2))^2-
\left(I\left((f+g)|f-g|\right)\right)^2+2I(f^2)I(g^2)\right].
$$
Тогда разница между правой и левой частями в неравенстве \eqref{10}  равна

$$
\mbox{RHS - LHS =} I(f^2)I(g^2)-
$$
$$
-\frac{1}{4}\left[(I(f^2))^2+(I(g^2))^2-
\left(I\left((f+g)|f-g|\right)\right)^2+2I(f^2)I(g^2)\right]=
$$
$$
=\frac{1}{4}\left[I((f+g)|f-g|)\right]^2 - \frac{1}{4}\left[I(f^2)-I(g^2)\right]^2=
$$
$$
=\frac{1}{4}\cdot\left[I((f+g)|f-g|)+I((f+g)(f-g))\right]\cdot
$$
$$
\cdot\left[I((f+g)|f-g|)-I((f+g)(f-g))\right]=
$$
$$
=\frac{1}{4}\cdot I\left((f+g)((f-g)+|f-g|)\right)\cdot I\left((f+g)((f-g)-|f-g|)\right)=
$$
$$
=\int_{E(f\ge g)} (f+g)(f-g)\,dx\cdot
\int_{E(f\le g)} (f+g)(g-f)\,dx=
$$
$$
=\int_{E(f\ge g)} (f^2-g^2)\,dx\cdot
\int_{E(f\le g)} (g^2-f^2)\,dx.
$$

Очевидно, что последнее выражение неотрицательно, таким образом неравенство \eqref{10} доказано.

\hfill $\Box$

Проведённые преобразования позволяют получить важные следствия.

\begin{corollary} \label{cor1}
Пусть $f(x), g(x)$  являются функциями любого знака, необязательно положительными. Тогда справедливо следующее тождество для разности между правой и левой частями в \eqref{10}

$$
\mbox{RHS  - LHS  =}
$$
$$
\int^b_af^2(x)dx\cdot \int^b_ag^2(x)\,dx - \int^b_a[\max(f,g)]^2dx\cdot\int^b_a[\min(f,g)]^2\,dx=\label{22}
$$
$$
\nonumber=\int_{E(f\ge g)} (f^2-g^2)\,dx\cdot
\int_{E(f\le g)} (g^2-f^2)\,dx.
$$
\end{corollary}

Подчеркнём, что последнее тождество справедливо для функций любого знака. Из него следуют некоторые обобщения неравенств Коши--Буняковского для необязательно неотрицательных функций. Приведём соответствующий результат.

Следующее замечание принадлежит Александру Артюшину.

\begin{corollary} Пусть для функций  $f(x), g(x)$ выполняется условие  $f(x) + g(x)\geq 0$ .  Тогда цепочка неравенств \eqref{10} также выполняется.
\end{corollary}

Из приведённого выше тождества также следует, что при противоположном условии $f(x) + g(x)\leq 0$
неравенства \eqref{10} выполняются в другую сторону.

Во всех рассмотренных случаях можно также явно указать, когда в неравенствах достигаются знаки равенства.

Отметим, что в работах автора  \cite{Sit1,Sit3,Sit4,Sit2} получены дальнейшие обобщения рассматриваемых задач. В том числе в  цитированных работах рассмотрено другое семейство обобщений неравенства Коши--Буняковского в терминах логарифмических производных, а также нетривиальный вопрос о сравнении различных уточнений между собой.

\section{Некоторые приложения}

Рассмотренные уточнения интегрального неравенства Коши--Буняковского могут быть перенесены на более общие случаи: бесконечные промежутки интегрирования, функции, интегрируемые по Риману или Лебегу, функции нескольких переменных и интегралы по многомерным областям. Следует заметить, что аккуратные доказательства в  этих случаях требуют зачастую достаточно кропотливых рассуждений.

Результаты также достаточно очевидным образом могут быть перенесены на уточнения дискретной и интегральной версии неравенства Минковского.

\subsection{  Принцип неопределённости для дискретного преобразования Фурье}

Рассмотрим дискретное преобразование Фурье (ДПФ), которое определяется на $n$--мерных векторах по формулам
$$
b_j=\frac{1}{\sqrt n}\sum_{k=0}^{n-1}a_k w^{-jk},
a_m=\frac{1}{\sqrt n}\sum_{k=0}^{n-1}b_k w^{mk},
$$
где $w=\exp(\frac{2\pi i}{n})$--первообразный корень из единицы.

Сейчас не будем агитировать, что это одно из самых важных понятий в математике, а приведём такую форму принципа неопределённости, которая доказана в \cite{DoHu}. Доказательство там не очень длинное, но я бы мягко говоря не назвал его понятным. Случай равенства указан один из возможных, но не все. Дадим другое доказательство, полученное Л.А.Мининым, которое получается в пару строк из неравенства Коши -- Буняковского.

\begin{theorem} (Л.А.Минин).  Пусть $A,B$--число ненулевых компонент в соответствующих векторах, $n$--их размерность. Тогда справедливо соотношение неопределённости
$AB\ge n$, причём равенство достигается тогда и только тогда, когда в одном из векторов все ненулевые компоненты равны.
\end{theorem}

Доказательство. Получаем, последовательно применяя неравенства треугольника и Коши--Буняковского
$$
|b_j|\le \frac{1}{\sqrt n}\sum_{a_k\neq0}|a_k|\le \frac{1}{\sqrt n}\sqrt A
\sqrt{\sum_{a_k\neq0}|a_k|^2}\le \sqrt{\frac{A}{n}}\|a\|_2.
$$
Далее возводим в квадрат, складываем ненулевые компоненты, применяем равенство Парсеваля--Стеклова--- вот и всё!!! Условия равенства  следуют из соответствующих условий для неравенства Коши--Буняковского.

Смысл этого утверждения в следующем: последовательность и её ДПФ не могут иметь одновременно слишком много нулей. 

Теперь можно выписывать различные обобщения приведённого неравенства по разработанной методике.

\subsection{Неравенства для гамма--функций и неполных гамма--функций}

 Рассмотрим приложение доказанного уточнения интегрального неравенства Коши--Буняковского (\ref{10}) для вывода оценок для гамма--функций и неполных гамма--функций. Напомним определения этих функций:
\begin{equation*}
\Gamma(a)=\int_0^{\infty} t^{a-1}e^{-t}\,dt,
\Gamma(a,x)=\int_x^{\infty} t^{a-1}e^{-t}\,dt, \gamma(a,x)=\int_0^x t^{a-1}e^{-t}\,dt,
\end{equation*}

Подставим в полученное уточнение неравенства Коши--Буняковского (\ref{10}) величины:
\begin{equation*}
f(x)=x^{\frac{a+1}{2}} e^{-\frac{x}{2}}, g(x)=x^{\frac{a-1}{2}} e^{-\frac{x}{2}},a>0.
\end{equation*}
Нетрудно проверить, что в этом случае
\begin{equation*}
\max(f,g)=
\begin{cases}
f(x), x\geq 1 , & \\
g(x), x\leq 1 , &
\end{cases}
\min(f,g)=
\begin{cases}
g(x), x\geq 1 , & \\
f(x), x\leq 1 . &
\end{cases}
\end{equation*}
Тогда указанное неравенство (\ref{10}) принимает такой вид
\begin{gather}\label{gamma1}
\Gamma^2(a+1)\leq
\left(\gamma(a,1)+\Gamma(a+2,1)\right)\cdot\left(\gamma(a+2,1)+\Gamma(a,1)\right)\leq
\Gamma(a+2)\cdot\Gamma(a).
\end{gather}

В результате нами доказано усиление известного свойства логарифмической выпуклости гамма--функции Эйлера. Компьютерное исследование показывает, что неравенство (\ref{gamma1}) является довольно точным, особенно для больших значений $a$, как и сама первоначальная оценка для логарифмической выпуклости. Неравенство (\ref{gamma1}) --- это также представитель известного семейства так называемых неравенств Турана, выражающих логарифмическую выпуклость  (или вогнутость) различных специальных функций по параметрам, см., например, \cite{KS1,Tu1,Tu2}.

Чтобы продемонстрировать точность полученного неравенства \eqref{gamma1} численно, преобразуем его, поделив на крайнюю правую часть. В результате получим с учётом известной формулы для гамма функции
$\Gamma(t+1)=t \Gamma(t)$ следующие неравенства
\begin{equation}\label{gamma2}
l(a)= \frac{a}{a+1} \leq
\frac{\left(\gamma(a,1)+\Gamma(a+2,1)\right)\cdot \left((\gamma(a+2,1)+\Gamma(a,1)\right)}{\Gamma(a+2)\Gamma(a)}=g(a) \leq 1.
\end{equation}
В этих неравенствах \eqref{gamma2} нижняя функция $l(a)$ не представляет интереса, интерес представляет скорость и точность приближения верхней функции $g(a)$ к единице по мере роста параметра $a$. Результаты отражены в приведённой  ниже таблице. В начале для небольших значений параметра $a$ вычисления велись с 10 разрядами после запятой, затем число разрядов постепенно увеличивалось, что видно из таблицы. Элемент таблицы вида $0,9(19)...$ означает число, в котором после нуля и десятичной запятой стоят 19 цифр 9 подряд.

\begin{table}[h!]
\begin{center}
\begin{tabular}{|l|c|} \hline
a & g(a)  \\ \hline
3 & 0,9665383895...  \\ \hline
5 & 0,9988760963...  \\ \hline
7 & 0,9999801314... \\ \hline
 10 &  0,9999999802... \\ \hline
20 & 0,9(19)...  \\ \hline
 30 & 0,9(34)... \\ \hline
50 & 0,9(64)...  \\ \hline
 100 & 0,9(159)... \\ \hline
\end{tabular}
\end{center}
\caption{Значения функции $g(a)$.}
\end{table}

Из приведённой таблицы видна очень хорошая  точность неравенства \eqref{gamma2}.

\subsection{Итерационные процедуры уточнения оценок для специальных функций}

В работах автора \cite{Sit1,Sit3,Sit4,Sit2} доказаны общие теоремы,  позволяющие каждое из
полученных ранее обобщений неравенства Коши-Буняковского развить в
некоторую итерационную процедуру,  генерирующую бесконечное число
новых последовательных обобщений того же неравенства, каждое
следующее из которых в определённом смысле точнее предыдущего.

При этом  основную роль начинают играть именно левые
очевидные (но не тривиальные!) части наших обобщений \eqref{D5}--\eqref{D55} (в дискретном случае) и \eqref{6}--\eqref{7} (в интегральном случае).
Основная идея проста: раз неравенство Коши--Буняковского позволяет
уточнять само себя, то этот приём можно просто повторять сколько
угодно раз. Как показано ниже, на этом пути удаётся, например,
получать эффективные оценки специальных функций, заданных своими
интегральными представлениями. Таким образом, как в хорошем
хозяйстве ничего не остаётся без полезного применения, так
оказываются востребованными и левые части обобщений \eqref{D5}--\eqref{D55}, \eqref{6}--\eqref{7}.

Процитируем одну из соответствующих теорем  из \cite{Sit1,Sit3,Sit4,Sit2} в наиболее простом варианте.

\begin{theorem}\label{proc}
 Определим итерационную процедуру уточнения
неравенства Коши - Буняковского по формулам:
\begin{equation}
\left\{
\begin{array}{ll}
U_0=f, V_0=g,\smallskip\nonumber\\
U_{n+1}=M_{\alpha}(U_n, V_n), V_{n+1}=M_{-\alpha}(U_n, V_n),\smallskip\nonumber\\
L_n=\int^b_a V_n^2\,dx, G_n=\int^b_aU_n^2\,dx,\smallskip\nonumber\\
A_n=L_nG_n.
\end{array}
\right.
\end{equation}
Тогда справедливы соотношения
\begin{equation}\label{it1}
S^2=\left(\int^b_afg\,dx\right)^2\le A_n=L_nG_n\le A_{n-1}\le
\int^b_af^2\,dx\int^b_ag^2\,dx,
\end{equation}
\begin{equation}\label{it2}
L_{n-1}\le L_n\le  S  \le G_n \le G_{n-1},
\end{equation}
\begin{equation}\label{it3}
S^2=\lim_{n\rightarrow\infty}A_n,
\end{equation}
где $M_\alpha$ - степенное среднее, $\alpha>0$.
\end{theorem}

Из неравенств \eqref{it1}--\eqref{it2} последней теоремы следует, что мы получаем
всё более и более точные двусторонние оценки левой части неравенства
Коши~- Буняковского, а из предельного соотношения \eqref{it3} следует, что
мы можем приближаться к левой части этого неравенства сколь угодно близко с обеих
сторон. Таким образом, если указанная левая часть --- это
интегральное представление некоторой специальной функции, то будут
получены сколь угодно точные двусторонние оценки для этой функции.

Аналогичные итерационные процедуры могут быть построены с
использованием средних Радо или с попеременным использованием
различных типов средних за счет  свободы их выбора. Ограничением
является необходимость явно вычислять интегралы из \eqref{it1}.

 На основании полученных результатов были разработаны итерационные процедуры для последовательного уточнения интегральных неравенств, которые могут быть приложены к различным специальным функциям. В качестве примера рассмотрим применение упомянутых
итерационных процедур к получению неравенств для полных
эллиптических интегралов Лежандра первого рода. Эта
гипергеометрическая функция имеет интегральное представление при
$0\le x< 1$ следующего вида:
\begin{equation}\label{Leg1}
K(x)=\int_0^1\frac{dt}{\sqrt{(1-t^2)(1-x^2t^2)}}
\end{equation}
и выражается через гипергеометрическую функцию Гаусса по формуле
\begin{equation*}
K(x)=\frac{\pi}{2}\cdot
{}_2F_1\left(1/2,1/2;1; x^2\right)
\end{equation*}
При  $x\rightarrow 1$ эллиптический интеграл Лежандра имеет особенность
с асимптотикой
\begin{equation}\label{Leg2}
K(x)=\frac{1}{2}\ln\frac{1}{1-x}+\frac{3}{2}\ln 2+o(1-x).
\end{equation}
Причём рост к бесконечности эллиптического интеграла Лежандра при $x\rightarrow 1-0$ очень медленный, например $K(0.9999999999)\approx 12.8992$.

Отметим, что разработаны многочисленные методы доказательства различных неравенств для эллиптических интегралов, многие из них достаточно сложные и трудоёмкие. Однако, стоит отметить, что достаточны точные неравенства можно получить в пару строк элементарными методами, причём даже с сохранением правильной асимптотики вблизи особенности. Например, используя неравенство между средними арифметическим и геометрическим, получаем для полного эллиптического интеграла Лежандра первого рода из определения \eqref{Leg1}

$$
\sqrt{(1-x^2)(1-k^2x^2)} \le \frac{1-x^2+1-k^2x^2}{2}=\frac{2-(1+k^2)x^2}{2}\ ,
$$
следовательно, 
$$
\frac{1}{\sqrt{(1-x^2)(1-k^2x^2)}} \ge \frac{2}{2-(1+k^2)x^2.}
$$

Интегрируя полученное неравенство, получаем такую оценку эллиптического интеграла снизу

$$
K(x) \ge \int_0^1 \frac{2 \,dx}{2-(1+k^2)x^2} =
\sqrt{
\frac{1}{2(1+k^2)}
}
\ \ln{
\frac
{
1+\sqrt{\frac{1+k^2}{2}}
}
{
1-\sqrt{\frac{1+k^2}{2}}}
}
$$

Последнее неравенство элементарно и получено элементарным методом, при этом вычисления демонстрирует, что оно достаточно точное, в том числе вблизи особенности $x\approx 1$.

Вернёмся к приведённым выше способам последовательного уточнения неравенств. 
Для эффективного проведения  итерационной
процедуры теоремы \ref{proc} ключевую роль играет правильное "расчленение" \
интеграла из \eqref{Leg1} в скалярное произведение. Осуществим следующий
выбор:
\begin{eqnarray*}
K^2(x)=\left(\int_0^1f(t)g(t)dt\right)^2,\\
f=\frac{1}{(1+t)^{1/2}(xt^2-(x+1)t+1)^{1/4}},\\
g=\frac{1}{(1+xt)^{1/2}(xt^2-(x+1)t+1)^{1/4}}.
\end{eqnarray*}
Он обусловлен необходимостью: во-первых, предусмотреть сходимость
интегралов на всех шагах итерационной процедуры; во-вторых,
выполнить равенство $f=g$ при $x=1$, которое позволяет надеяться на
повышенную точность оценок при $x\rightarrow 1$; в-третьих,  явно
вычислять получаемые интегралы.

Мы приведем
список вычисленных оценок итерационной процедуры \ref{proc}, при условии, что выполнены три первых шага этой процедуры, то есть вычислены нижние оценки $L_0,L_1,L_2$ и верхние оценки $G_0,G_1,G_2$ для эллиптического интеграла Лежандра $K(x)$.   При этом нетривиальные оценки получаются уже из самого неравенства Коши - Буняковского на начальном шаге  за счет удачного выбора пары  функций $f,g$.

\begin{theorem}
Справедливы  следующие неравенства для
полного эллиптического интеграла Лежандра  первого рода  $K(x)$,  следующие из итерационной процедуры \ref{proc}:
\begin{equation}\label{Leg3}
L_0(x) \leq L_1(x) \leq L_2(x) \ldots \leq K(x) \leq \ldots G_2(x) \leq G_1(x) \leq G_0(x).
\end{equation}
\end{theorem}
При этом указанные оценки вычисляются по формулам
\begin{gather*}
L_0(x)=\frac{1}{\sqrt{2(x+1)}}\ln\left(\frac{2\sqrt{2(x+1)}+x+3}{1-x}\right),
\\
G_0(x)=\frac{1}{\sqrt{(x+1)2x}}\ln\left(\frac{2\sqrt{2x(x+1)}+3x+1}{1-x}\right),
\\
L_1(x)=\frac{2}{\sqrt{(x+3)(3x+1)}}\ln\left(\frac{\sqrt{(x+3)(3x+1)}+2x+2}{1-x}\right),
\\
G_1(x)=\frac{1}{2}\left[\frac{1}{\sqrt{2x(x+1)}}
\ln\left(\frac{2\sqrt{2x(x+1)}+3x+1}{1-x}\right)\right]+\\
+\frac{1}{\sqrt{2(x+1)}}\ln\left(\frac{2\sqrt{2(x+1)}+x+3}{1-x}\right),
\\
L_2(x)=\frac{1}{\sqrt{2x+2}}\ln\left(\frac{2\sqrt{2x+2}+x+3}{1-x}\right)+ \\
+\frac{1}{\sqrt{2x(x+1)}}\ln\left(\frac{2\sqrt{2x(x+1)}+3x+1}{1-x}\right)+ \\
+\frac{1}{\sqrt{(x+3)(3x+1)}}\ln\left(\frac{\sqrt{(x+3)(3x+1)}+2x+2}{1-x}\right).\\
G_2(x)=\frac{\left(\frac{\sqrt{5}+1}{2\sqrt{5}}\right)} { \sqrt{
(x+4+\sqrt{5}) \left( \left(\frac{7+\sqrt{5}}{2}
\right)x+\frac{3+\sqrt{5}}{2}\right) }
} \times\\
\times\ln\left(\frac{2\sqrt{(x+4+\sqrt{5})\left(\left(\frac{7+\sqrt{5}}
{2}\right)x+\frac{3+\sqrt{5}}{2}\right)}
+\left(\frac{9+\sqrt{5}}{2}\right)x+\frac{11+3\sqrt{5}}{2}}{\left(\frac{5+\sqrt{5}}
{2}\right)(1-x)}\right)+\\
+\frac{ \left( \frac{\sqrt{5}-1}{2\sqrt{5}} \right) } { \sqrt{
(x+4-\sqrt{5}) \left( \left( \frac{7-\sqrt{5}} {2} \right)
x+\frac{3-\sqrt{5}}{2} \right) } }\times\\
\times\ln\left(\frac{2\sqrt{(x+4-\sqrt{5})\left(\left(\frac{7-\sqrt{5}}
{2}\right)x+\frac{3-\sqrt{5}}{2}\right)}
+\left(\frac{9-\sqrt{5}}{2}\right)x+\frac{11-3\sqrt{5}}{2}}{\left(\frac{5-\sqrt{5}}
{2}\right)(1-x)}\right),
\end{gather*}

Отметим, что первоначальные двусторонние оценки с величинами $L_0(x), G_0(x)$ получаются применением обычного неравенства Коши--Буняковского. 

Расчеты показывают высокую точность полученных двусторонних
неравенств \eqref{Leg3}. Можно продолжить применённую итерационную процедуру и получить любое число всё более точных двусторонних неравенств для полных эллиптических интегралов Лежандра первого рода. Тематике неравенств для эллиптических интегралов Лежандра всех трёх родов посвящено достаточно большое число работ, см., например, \cite{KS2,KS3}.

\subsection{Уточнения неравенства Ацеля в пространстве Лоренца}

Сравнительно малоизвестным является тот факт, что в пространствах
Лоренца со знаконеопределённой метрикой (или Понтрягина, см. \cite{BB}) также выполняется неравенство Коши--Буняковского, но в
обратную сторону.  Данное
неравенство является частным случаем известных неравенств Ацеля. Как указано в \cite{BB} именно неравенство Ацеля позволяет математически строго обосновать так называемый "парадокс близнецов"\ в специальной теории относительности Эйнштейна.

Для этого случая получена

\begin{theorem}
Пусть $A(x, y)$ есть средняя часть
произвольного обобщения дискретного неравенства Коши--Буняковского
вида
$$
\left(\sum_{k=1}^nx_k\cdot y_k\right)^2\le A(x,
y)\le\left(\sum_{k=1}^nx_k^2\right)\cdot\left(\sum_{k=1}^ny_k^2\right),
$$
и выполнены условия
$$
x_0^2-\sum_{k=1}^nx_k^2\geq 0,\  y_0^2-\sum_{k=1}^ny_k^2\geq 0.
$$
Тогда справедливо следующее уточнение дискретного неравенства
Коши--Буняковского в пространствах Лоренца:
\begin{gather*}
\left(x_0y_0-\sum_{k=1}^nx_ky_k\right)^2\geq \left(x_0y_0-\sqrt{A(x,
y)}\right)^2 \ge\\
\ge\left(x_0^2-\sum_{k=1}^nx_k^2\right)\left(y_0^2-\sum_{k=1}^ny_k^2\right).
\end{gather*}
\end{theorem}

 Аналогично обычному случаю из последнего неравенства выводится
обобщение неравенства Минковского для пространств Лоренца.

\subsection{Уточнения неравенства Коши--Буняковского для $q$ интегралов Джексона }
Следующий пример относится к области коммутативного $q$--Анализа. Это такой раздел Анализа, в котором есть свои понятия $q$--производной, $q$--интеграла и соответствующие обобщения специальных функций, которые при $q\to 1$ переходят в обычные определения, см. \cite{GR}. Так, например, $q$--интеграл Джексона вводится по формуле
$$
\int_0^1 f(t)\,d_q t=(1-q)\sum_{k=0}^\infty f(q^k) q^k.
$$
Этот интеграл на самом деле не интеграл, а ряд.
Для него наши результаты позволяют получить такое обобщение неравенства Коши--Буняковского в терминах $q$ - интеграла Джексона.

\begin{theorem}
 Пусть $M$ - произвольное абстрактное среднее‚ $M^*$ -
сопряженное к нему. Тогда справедливо обобщение неравенства Коши--Буняковского в терминах $q$ - интеграла Джексона:
\begin{align*}
\left(\int_0^1f(t)g(t)d_qt\right)^2 &\le \left(\int_0^1(M(f(t),
g(t)))^2d_qt\right)\cdot\\
\cdot\left(\int_0^1(M^*(f(t),
g(t)))^2d_qt\right) &\le
\left(\int_0^1f^2(t)d_qt\right)^2\left(\int_0^1g^2(t)d_qt\right)^2.
\end{align*}
\end{theorem}

\subsection{Задача об обобщениях неравенств Коши--Буняковского и Ацеля для знаконеопределённой формы от четырёх переменных в пространстве Понтрягина}

 Рассмотрим задачу о возможных обобщениях дискретных неравенств Коши--Буняковского и Ацеля для знаконеопределённой формы от четырёх переменных в пространстве Понтрягина (или Лоренца).
Речь идёт о квадратичных неопределённых формах с сигнатурой знаков $(+,\ +,\ -,\ -\ )$ в пространствах Понтрягина-Лоренца.

\textbf{Задача:} описать области переменных в $\mathbb{R}^4 \times \mathbb{R}^4$ (или $\mathbb{C}^4 \times \mathbb{C}^4$), когда нижеприведённая форма положительна, отрицательна, равна нулю.\\
Это одновременное обобщение неравенств Коши--Буняковского и Ацеля.
\begin{gather*}
 (x_1y_1+x_2y_2-x_3y_3-x_4y_4)^2 -
  (x_1^2+x_2^2-x_3^2-x_4^2)(y_1^2+y_2^2-y_3^2-y_4^2)=  \\
  =+ y_2^{2} x_4^{2}+ y_1^{2} x_3^{2}+ x_2^{2} y_3^{2}+ x_1^{2}
     y_3^{2}+
       y_1^{2} x_4^{2}+ x_1^{2} y_4^{2}+ y_4^{2} x_2^{2}+ y_2^{2}
     x_3^{2}-\\
    - y_4^{2} x_3^{2}- y_3^{2} x_4^{2}- x_1^{2} y_2^{2} - y_1^{2}
    x_2^{2}+\\
                                                              +2  y_4 x_3 y_3 x_4+2  x_1 y_1 x_2 y_2-\\
                                                            -2  x_1 y_1 x_3 y_3-2  x_2 y_2 x_3 y_3-2  x_1 y_4 y_1 x_4-2  y_4
                                                                  x_2
                                                              y_2 x_4=
  \\
    ={( y_3 x_1- x_3 y_1)}^{2}+{( x_4 y_1- x_1 y_4)}^{2}+{( y_3 x_2-
       y_2
       x_3)}^{2}+{( y_2 x_4- y_4 x_2)}^{2}-\\
                                                            -{( y_2 x_1- y_1 x_2)}^{2}-{( y_3 x_4- x_3 y_4)}^{2}
\end{gather*}

Решение этой задачи автору неизвестно. Для возможных подходов к решению  может оказаться полезным применённое в последнем равенстве тождество Лагранжа для комплексных величин (замечание Владимира Кисиля). К сожалению, кажущийся очевидным путь решения этой задачи на основании приведённого тождества Лагранжа сразу не приводит цели, так как приравненные нулю формы под знаками квадратов оказываются зависимыми.

Было бы интересным установить геометрический смысл действительных или комплексных многомерных поверхностей, появляющихся в формулировке этой задачи. Разумеется, можно рассматривать аналогичные задачи для знаконеопределённых форм с произвольной сигнатурой в пространствах Понтрягина. Решения этих задач представляются автору чрезвычайно трудными, возможно, хоть какое-то продвижение в них окажется возможным с использованием компьютерного моделирования.

\subsection{ Принцип переноса }

Отметим, что существует разработанная А. Г. Кусраевым техника перенесения результатов для числовых неравенств на гораздо
более общий случай векторных решёток и билинейных операторов
в них.  При этом основным инструментом становится приближение выпуклых функций с помощью прямых, как в классических теоремах Минковского и Хёрмандера из выпуклого анализа, см. \cite{Kus}. Приведём цитату из этой работы.

"In a private discussion Professor S.M. Sitnik hypothesized that there must exist some general principle allowing to produce automatically new inequalities for bilinear operators, provided that the corresponding ones hold true for bilinear forms.
The aim of this paper is to present a transfer principle which enables us to transform inequalities with semi-inner products to inequalities containing positive semidefinite symmetric bilinear operators with values in a vector lattice."

Таким образом, в силу принципа переноса можно единым универсальным способом распространять обычные неравенства для функций на случай более общих абстрактных конструкций.

\section{Заключение и благодарности}

Таким образом, мы постарались показать, что тематика уточнений классических неравенств является живой и развивающейся областью современной математики. К этому направлению, безусловно, относятся и уточнения интегрального и дискретного неравенств Коши--Буняковского. Часть результатов этого направления стали классическими и получили различные приложения, например, как теорема CDE --- Дейкина--Карлица--Элиезера. Получаемые как обобщения классических неравенств на этом пути оценки находят многочисленные приложения в теории дифференциальных уравнений, операторах преобразования, теории интегральных преобразований и специальных функций, а также во многих других разделах современной математики. К сожалению, многие достижения  в этой области мало известны широкому математическому сообществу. Целью данной статьи было хотя бы частично восполнить указанный пробел.

\textbf{Благодарности.} Автор благодарит за плодотворные обсуждения  вопросов этой работы и различную помощь в различные моменты времени Ларису Ивановну Брылевскую (Санкт--Петербург), Анатолия Георгиевича Кусраева (Владикавказ), Анатолия Кореновского (Одесса), Александра Артюшина (Новосибирск), Владимира Кисиля (Лидс).


\bigskip

\begin{thebibliography}{777}
	
	\bibitem{BB} %
Э. Беккенбах, Р. Беллман. {\it  Неравенства},  М.: Мир, 1965.

\bibitem{Jini} %
 К. Джини. {\it  Средние величины}, Москва: Статистика, 1970.

\bibitem{Kal} %
 С.И. Калинин. {\it  Средние величины степенного типа. Неравенства Коши и Ки
Фана},  Киров, 2002.

\bibitem{T1} %
 В.В. Катрахов, С.М. Ситник. {\it Метод операторов преобразования и краевые задачи для сингулярных эллиптических уравнений}, Современная математика. Фундаментальные направления, \textbf{64}:2 (2018), 211--426.

\bibitem{Met} %
Л.А. Ковалева, С.М. Ситник, О.В. Чернова, Э.Л. Шишкина. {\it
Математический анализ. Функции одной переменной.
Неравенства},  учебно-методическое пособие. Белгород: ИД «БелГУ», НИУ «БелГУ», 2021.

\bibitem{Mas1} В.~П. Маслов. {\it Нелинейное финансовое осреднение, эволюционный процесс и законы эконофизики}, Теория вероятностей и приложения, \textbf{49}:2 (2004), 269–-296.
\bibitem{Mas2} В.~П. Маслов. {\it О нелинейности осреднений в финансовой математике}, Математические заметки, \textbf{74}:6 (2003),944–-947.
\bibitem{Mas3} В.~П. Маслов. {\it Аксиомы нелинейного осреднения в финансовой математике и динамика курса акций}, Теория вероятностей и приложения, \textbf{48}:4 (2003), 800–-810.

\bibitem{MS} %
Л.А. Минин, С.М. Ситник.  {\it О неравенствах для тета-функций Якоби},
Чернозёмный альманах научных исследований. Серия "Фундаментальная математика".
Номер, посвящённый 85-летию Ивана Александровича Киприянова, \textbf{1}:8 (2009),  234--311.

\bibitem{PS} %
А.Б. Певный, С.М. Ситник {\it
Обобщения неравенств M. Г. Крейна, Е. А. Горина и Ю. В. Линника для положительно
определенных функций на многоточечный случай},
Сибирские электронные математические известия,  \textbf{16} (2019), 263--270.

\bibitem{SKM} %
С.~Г. Самко, А.~А. Килбас, О.~И. Маричев. {\it Интегралы и производные дробного порядка и некоторые их приложения}, Минск: Наука и техника, 1987.

\bibitem{Sit1} %
С.~М. Ситник.  {\it Уточнения и обобщения классических неравенств}.  Математический форум. Исследования по математическому анализу,  2009, 221--266.

\bibitem{Sit3} %
С.~М. Ситник.  {\it Обобщения неравенств Коши-Буняковского методом средних значений и их приложения},
 Черноз\"{е}мный альманах научных исследований. Серия "Фундаментальная математика", \textbf{1}:1 (2005), 3--42.

\bibitem{T2} %
С.М. Ситник, Э.Л. Шишкина. {\it Метод операторов преобразования для дифференциальных уравнений с операторами Бесселя}, М.: Физматлит, 2019.

\bibitem{WW} %
Э.~Т. Уиттекер, Дж.~Н. Ватсон. {\it Курс современного Анализа.  Часть вторая. Трансцендентные функции}, М.:~ГИФМЛ, 1963.

\bibitem{HLP} %
Г. Харди, Дж. Литтлвуд, Г. Полиа. {\it Неравенства}, М.: ИЛ,  1948.

\bibitem{Sit4} %
 P. Agarwal, A. Korenovskii, S. Sitnik. {\it A Generalization of Cauchy--Bunyakovsky Integral Inequality Via Means with Max and Min Values},
Chapter 18.\ In:  Trends in Mathematics. Advances in Mathematical Inequalities and Applications.
Eds.:  P. Agarwal, S.S. Dragomir,  M. Jleli,  B. Samet.
Birkhauser Basel, Springer Nature Singapore,
2018, 333--349.

\bibitem{BMV} P.~S. Bullen, D.~S. Mitrinovi\'{c}, P.~M. Vasi\'{c}. {\it Means and Their Inequalities}, D.Reidel Publishing Company: Dordrecht, 1988.

\bibitem{Bec} %
W. Beckner. {\it Inequalities in Fourier analysis}, Annals of Mathematics, \textbf{102}: 6 (1975), 159--182.

\bibitem{Bun} %
V. Buniakowski. {\it Sur quelques in\'{e}galit\'{e}s concernant les
int\'{e}grales ordinaires et les int\'{e}grales aux diff\'{e}rences
finies},  M\'{e}moires de l' Acad. de St. -
P\'{e}tersbourg (VII), \textbf{1}: 9  (1859).

\bibitem{Fum} %
 H. Fumio. {\it  Means of Hilbert Space Operators},
 Springer Lecture Notes in Mathematics, № 1820,  2000.

\bibitem{Dr1} %
S. Dragomir. {\it A Survey on Cauchy--Buniakowsky--Schwartz Type
Discrete Inequalities},  RGMIA monographs, 2003.

\bibitem{Har} %
P. Hartman. {\it Convex functions and mean value
inequalities},  Duke Math. J.,  \textbf{39}: 2 (1972),  351 -- 360.

\bibitem{GR} %
G. Gasper, M. Rahman. {\it Basic Hypergeometric Series}, 2nd edition, Cambridge University Press, 2004.

\bibitem{KS1} %
D. Karp, S.M. Sitnik. {\it Log-convexity and log-concavity of hypergeometric-like functions},
Journal of Mathematical Analysis and Applications. Elsevier, Amsterdam, \textbf{364} : 2 (2010), 384--394.

\bibitem{KS2} %
D. Karp, S.M. Sitnik. {\it
Asymptotic approximations for the first incomplete elliptic integral near logarithmic singularity}, Journal of Computational and Applied Mathematics. Elsevier, Amsterdam. \textbf{205} : 1 (2007), 186--206.

\bibitem{KS3} %
D. Karp, A. Savenkova, S.M. Sitnik. {\it
Series expansions for the third incomplete elliptic integral via partial fraction decompositions},
Journal of Computational and Applied Mathematics. Elsevier, Amsterdam.
\textbf{207} : 2 (2007) 331--337.
Proceedings of The Conference in Honour of Dr. Nico Temme on the Occasion of his 65th birthday.

\bibitem{T4} %
Ed. V.V. Kravchenko, S.M. Sitnik. {\it Transmutation Operators and Applications},  In the Series: Trends in Mathematics.  Springer, Birkhauser, 2020.

\bibitem{MO2} %
A. Marshall, I. Olkin, B. Arnold. {\it Inequalities: Theory of Majorization and Its Applications.  Second Edition}, Springer, 2011.

\bibitem{Tu1} %
Mehrez K.,  Sitnik S. M. {\it Turan Type Inequalities for Classical and Generalized Mittag-Leffler Functions},
Analysis Mathematica, \textbf{44}: 4 (2018), 521--541.

\bibitem{MPF} %
D. Mitrinovi\'{c}, J. Pe\v{c}ari\'{c}, A. Fink. {\it Classical and new inequalities in analysis}, Kluwer, 1993.

\bibitem{MBV} %
D.S. Mitrinovi\'{c}.  {\it Means and their inequalities}, D.Reidel,
1988.

\bibitem{Rado} %
T. Rado.  {\it On convex functions}, Trans. Amer. Math. Soc.,  \textbf{37} (1935), 266 -- 285.

\bibitem{T3} %
E.L. Shishkina, S.M. Sitnik. {\it Transmutations, Singular and Fractional Differential Equations with Applications to Mathematical Physics},
    In the Series: Mathematics in Science and Engineering.
    Elsevier, Academic Press. 2020.

\bibitem{Sit2} %
S.~M. Sitnik. {\it Generalized Young and Cauchy--Bunyakowsky Inequalities with Applications: a survey}, arXiv:1012.3864,  2010,   51~p.

\bibitem{Tu2} %
S.M. Sitnik S.M., Kh. Mehrez. {\it
Proofs of some conjectures on monotonicity of ratios of Kummer, Gauss and generalized hypergeometric functions}, Analysis. (De Gruyter). \textbf{36}: 4 (2016), 263--268.

\bibitem{Ste} %
J. Steele. {\it The Cauchy-Schwarz Master Class: An Introduction to the Art of Mathematical Inequalities}, Cambridge University Press, 2004.

\bibitem{Toad} %
Gh. Toader. {\it Greek Means and the Arithmetic-Geometric Mean},  RGMIA monographs, 2005.

\bibitem{Kus} %
A.G. Kusraev. {\it A transfer principle for inequalities in vector lattices}, Journal of Mathematical Analysis and Applications. \textbf{374}: 1 (2011), 282--289.

\bibitem{DoHu} D.L. Donoho, X. Huo. {\it Uncertainty Principles and Ideal Atomic Decomposition}, Transactions on Information Theory. \textbf{47}: 7 (2001), 2845--2862.

\end{thebibliography}
\end{document}